%% file: main.tex
\documentclass[a4paper]{cas-sc}
\input{preamble}
\hypersetup{pdftitle={Kernel Compensation and Quasi-Newton Iteration for Nonlinear Non-Hermitian Eigenproblems in Photonic Crystals},pdfauthor={Chenhao Jin and Hehu Xie}}
\begin{document}
\let\WriteBookmarks\relax
\shorttitle{\mbox{}}
\shortauthors{C. Jin and H. Xie}
\title[mode=title]{Kernel Compensation and Quasi-Newton Iteration for Nonlinear Non-Hermitian Eigenproblems in Photonic Crystals}
\tnotemark[1]
\tnotetext[1]{This work was supported by National Science Challenge Project (TZ2025007), Natural Science Foundations of China (No. 12331015), National Key Laboratory of Computational Physics (No. 6142A05230501), Strategic Priority Research Program of the Chinese Academy of Sciences (XDB0640000, XDB0640300, XDB0620203, XDA0480504), National Key Research and Development Program of China (2023YFB3309104),  National Center for Mathematics and Interdisciplinary Science, CAS.}
\author[1]{Chenhao Jin}[orcid=0009-0001-5857-5400]
\ead{kanit@mail.ustc.edu.cn}
\affiliation[1]{organization={School of Mathematical Sciences, University of Science and Technology of China},city={Hefei},state={Anhui},postcode={230026},country={China}}
\author[2,3]{Hehu Xie}[orcid=0000-0002-1947-9256]
\cormark[1]
\ead{hhxie@lsec.cc.ac.cn}
\affiliation[2]{organization={SKLMS, Academy of Mathematics and Systems Science, Chinese Academy of Sciences},addressline={No.55, Zhongguancun Donglu},city={Beijing},postcode={100190},country={China}}
\affiliation[3]{organization={School of Mathematical Sciences, University of Chinese Academy of Sciences},city={Beijing},postcode={100049},country={China}}
\cortext[1]{Corresponding author}
\begin{abstract}
We extend the kernel-compensated eigensolver developed in our previous work from Hermitian to non-Hermitian nonlinear eigenproblems arising in dispersive photonic crystals. A mimetic finite-difference discretization preserves the discrete curl--divergence complex, and three-dimensional discrete Fourier transforms provide matrix-free realizations of the curl and divergence blocks. The isotropic inverse dielectric coefficient is stored through the material-edge indices, so its value and derivatives can be updated without assembling a Maxwell matrix. Kernel compensation is characterized by a fixed unitary block diagonalization, which remains exact for the nonlinear problem. At each outer step, a right Generalized Davidson iteration solves the non-Hermitian linearized eigenproblem. The exact derivative of a simple eigenvalue branch is expressed through its left and right eigenvectors. To avoid a second eigensolve, we replace the left eigenvector by the available normalized right eigenvector and obtain a right-Rayleigh quasi-Newton update. We prove that the exact left--right Newton iteration is locally quadratically convergent at simple roots, whereas the proposed iteration with exact inner eigenpairs has an explicit linear error factor and is at least quadratic when the derivative defect vanishes. A spectral regularity estimate bounds the gradient appearing in the derivative defect.
\end{abstract}
\maketitle
\section{Introduction}
\label{se:int}

Frequency-dependent permittivities turn the Maxwell eigenvalue problem in
photonic crystals into a nonlinear eigenvalue problem (NEP). Material
absorption leads to complex frequencies and non-Hermitian matrices
\cite{sigalas1994dispersive,effenberger2012absorbing}.
The resulting three-dimensional problems combine material nonlinearity
with a large curl kernel. Iterative frequency updates require repeated
solutions of large eigenproblems, so both the nonlinear iteration and
the treatment of the kernel affect the computational cost
\cite{guttel2017,huang2016newton}.

Several methods address the frequency dependence. Auxiliary-field and
linearization formulations replace the NEP by a larger linear problem
\cite{raman2010dispersive,effenberger2012absorbing}, while spectral
indicator methods treat frequency-dependent permittivities through
finite-element discretization \cite{xiao2021frequency}.
For general analytic NEPs, infinite Arnoldi uses an infinite-dimensional
linear operator \cite{jarlebring2012infinite}, NLEIGS uses rational
interpolation and structured linearization \cite{guttel2014nleigs}, and
contour integration recovers eigenvalues inside a prescribed contour
from resolvent moments \cite{beyn2012integral}.
Quasi-Newton analyses identify the effect of Jacobian approximations
\cite{jarlebring2018qn}; inexact Newton analyses relate inner accuracy
to local convergence \cite{dembo1982,szyld2013inexact}.
Our analysis concerns the particular derivative approximation used in
repeated frozen-frequency Maxwell eigensolves.

The Yee discretization preserves the discrete curl--divergence relation
\cite{yee1966}. FFT-based null-space-free methods exploit the curl
structure to remove zero eigenvalues \cite{huang2013eigendecomppc}.
Tiexiang Li and collaborators developed this approach in FAME, including CPU and
GPU implementations \cite{lyu2021fame}, and extended it to anisotropic
crystals in oblique coordinates \cite{tian2023oblique}.
Their structure-preserving complex-band method fixes a real frequency
and computes a complex wave-vector component through a palindromic
quadratic eigenvalue formulation \cite{huang2020complex}.
That formulation differs from the present NEP, which fixes a real Bloch
wave vector and seeks complex frequencies in a dispersive material.
For dispersive metallic crystals, Huang, Lin, and Mehrmann combined a
Newton-type iteration with nonequivalence deflation and FFT-based
null-space-free eigensolves \cite{huang2016newton}.
Their scalar Newton update uses a left eigenvector as well as the right eigenpair.

Our earlier work introduced kernel compensation with mimetic finite
differences \cite{jin2025kcpc}. A discrete divergence penalty treats the
curl kernel while retaining the original vector space.
We subsequently developed anisotropic discretizations and a
GPU-accelerated Hermitian solver \cite{paper2}, and combined kernel
compensation with fixed-point and Newton iterations for nonlinear
Hermitian problems \cite{paper3}.
In that setting, the normalized right eigenvector gives the exact
eigenvalue derivative, and the Hermitian min--max principle separates
the compensated spectrum.

The non-Hermitian extension requires different spectral and convergence
arguments. The min--max principle is unavailable for complex spectra,
and a simple eigenvalue derivative involves both left and right
eigenvectors \cite{greenbaum2020}.
We retain the mimetic DFT operators and use a right Generalized Davidson
eigensolver \cite{davidson1975,sleijpen1996} at each frozen frequency.
The quasi-Newton update uses the available normalized right eigenvector
in place of the left eigenvector and avoids an additional adjoint eigensolve. Its derivative defect is expressed through the
gradient of the material Rayleigh quotient and bounded using the
relative spectral gap of the Hermitian part.

The contributions of this work are as follows.
\begin{itemize}
\item A frequency-independent unitary transformation block-diagonalizes the
compensated non-Hermitian matrix, separating the discrete curl image,
divergence image, and harmonic spaces, and proves that kernel compensation
preserves the original nonlinear eigenpairs under an explicit penalty condition.
\item The spectral analysis of the preconditioner is extended to complex
material coefficients: the two-sided preconditioned matrix is normal, its
spectrum lies on a complex line segment, and its condition number admits a
material-dependent bound independent of the grid size.
\item The right-Rayleigh update reuses the computed right eigenvector.
For exact inner eigenpairs, we express its linear error factor in terms of
the gradient of the material Rayleigh quotient. A relative-gap estimate
for the Hermitian part bounds the derivative defect and gives a
sufficient condition for local contraction at a simple root.
We implement the method in complex double precision with matrix-free CUDA
operators. For the \texttt{SC\_CURV}/\texttt{SiC} tests with
$64\leq N\leq128$, the GPU implementation achieves speedups of
$69.6$--$73.0$ over the single-core, single-thread CPU implementation.
\end{itemize}

The paper is organized as follows. Section~\ref{se:disc} gives the mimetic
discretization and matrix-free operators. Section~\ref{se:kc} establishes
the kernel compensation formulation and preconditioning results.
Section~\ref{se:nep} presents the eigensolver, the quasi-Newton update,
and the local convergence analysis.
Section~\ref{se:num} reports complex band structures, convergence
histories, and CPU and GPU timings. Section~\ref{se:con} concludes the paper.

\section{Discretization}
\label{se:disc}

\subsection{Model Problem and Mimetic Finite Difference Scheme}
\label{sse:disc:mfd}

The model problem is defined over the primitive cell
\begin{align}
\label{display:dis:cell}
\Omega:=\left\{\vec x=\sum_{i=1}^3x_i\vec a_i:x_i\in[0,1],\ i=1,2,3\right\},
\end{align}
where $\{\vec a_i\}_{i=1}^3$ forms a basis in $\Rbb^3$. For a Bloch wave
vector $\wnk$, the original magnetic field satisfies
$$
\vec H^{\mathrm{qp}}(\vec x+\vec a_i)
=\exp(\imagunit\wnk\cdot\vec a_i)\vec H^{\mathrm{qp}}(\vec x).
$$
The phase transformation
$\vec H(\vec x):=\exp(-\imagunit\wnk\cdot\vec x)\vec H^{\mathrm{qp}}(\vec x)$
makes the transformed field periodic. The corresponding differential operator is
\begin{align}
\label{display:dis:shift}
\nabla_{\wnk}:=\nabla+\imagunit\wnk.
\end{align}
The magnetic-field formulation with this shifted operator is
\begin{align}
\label{display:dis:model}
\left\{
\begin{aligned}
\nabla_{\wnk}\times
\varepsilon^{-1}(\vec x,\omega)
\nabla_{\wnk}\times\vec H
&=\omega^2\vec H
\quad\text{in }\Omega,\\
\nabla_{\wnk}\cdot\vec H
&=0
\quad\text{in }\Omega,
\\
\vec H(\vec x+\vec a_i)
&=\vec H(\vec x),
\quad i=1,2,3.
\end{aligned}
\right.
\end{align}
The dielectric coefficient $\varepsilon$ satisfies the spatial periodicity condition: $\varepsilon(\vec x+\vec a_i,\omega)=\varepsilon(\vec x,\omega)$.  $\varepsilon^{-1}$ is piecewise defined in $\Omega$:
\begin{align}
\label{display:dis:diel}
\varepsilon^{-1}(\vec x,\omega)
=
\begin{cases}
\eta(\omega),&\vec x\in\Omega_1,\\
1,&\vec x\in\Omega\setminus\Omega_1,
\end{cases}
\end{align}
where $\Omega_1$ denotes the material subregion. The function $\eta$ is holomorphic in a pole-free neighborhood of each frequency under consideration. We consider frequencies in $\Omega_{\mathrm{sol}}\subset\Cbb$ such that
\begin{align}\label{display:dis:soldomain}
    \Omega_{\mathrm{sol}}\subset\{\omega:\ \omega\ne0,\ 0<\re\eta(\omega)<1\}.
\end{align}

We map the primitive cell to the unit cube and divide each coordinate direction into $N$ intervals of width $h=\frac{1}{N}$. As in our previous mimetic finite-difference discretizations \cite{jin2025kcpc,paper2,paper3}, the Yee scheme \cite{yee1966} places the three electric components at edge centers, the three magnetic components at face centers, and the divergence value at the cell center, as is illustrated in Table~\ref{tab:dis:dof}.
\begin{table}[pos={htb}]
    \centering
    \caption{Degrees of freedom of the discrete electromagnetic fields and divergence image on the Yee grid with $N$ intervals in each direction. The indices satisfy $i,j,k\in\{1,2,\ldots,N\}$.}
    \label{tab:dis:dof}
\normalfont\fontsize{11}{13.6}\selectfont
\renewcommand{\small}{\fontsize{11}{13.6}\selectfont}

\renewcommand{\arraystretch}{1.3}
\setlength{\extrarowheight}{2pt}

    \begin{tabular*}{\linewidth}{@{\extracolsep{\fill}}ccc@{}}
    \toprule
      Physical Quantity &  Functions & Approximate Vectors \\
      \midrule
      Electric Field & $\vec E=(E^1,E^2,E^3)$ & $\begin{pmatrix}
          E^1_{i-\frac{1}{2},j,k} & E^2_{i,j-\frac{1}{2},k} & E^3_{i,j,k-\frac{1}{2}}          
      \end{pmatrix}$ \\
      Magnetic Field & $\vec H=(H^1,H^2,H^3)$ & $\begin{pmatrix}
          H^1_{i,j-\frac{1}{2},k-\frac{1}{2}} & H^2_{i-\frac{1}{2},j,k-\frac{1}{2}} & H^3_{i-\frac{1}{2},j-\frac{1}{2},k}          
      \end{pmatrix}$ \\
      Divergence & $\psi$ & $\begin{pmatrix}
          \psi_{i-\frac{1}{2},j-\frac{1}{2},k-\frac{1}{2}}
      \end{pmatrix}$\\
      \bottomrule
    \end{tabular*}
\end{table}
\FloatBarrier
The symbol $\psi$ denotes the divergence. The discrete electric and magnetic fields each contain $3N^3$ degrees of freedom (DoFs), whereas the discrete divergence contains $N^3$ DoFs. The outer discrete curl maps edge DoFs to face DoFs; its conjugate transpose gives the inner curl. The discrete inverse dielectric matrix $\Mcal_\eta(\omega)$ maps edge DoFs to edge DoFs.

\subsection{Discrete Matrix Operators}
\label{sse:disc:mat}

The matrix representations of the discrete operators have circulant structures. We introduce two $N\times N$ circulant matrices:
\begin{align}
\label{display:dis:circulant}
\Dbf_1
:=\frac{1}{h}
\begin{pmatrix}
1&0&\cdots&0&-1\\
-1&1&\ddots&&0\\
0&\ddots&\ddots&\ddots&\vdots\\
\vdots&&-1&1&0\\
0&\cdots&0&-1&1
\end{pmatrix},
\qquad
\Dbf_0
:=\frac{1}{2}
\begin{pmatrix}
1&0&\cdots&0&1\\
1&1&\ddots&&0\\
0&\ddots&\ddots&\ddots&\vdots\\
\vdots&&1&1&0\\
0&\cdots&0&1&1
\end{pmatrix}.
\end{align}
The matrix $\Dbf_1$ approximates a first derivative between staggered locations, while $\Dbf_0$ averages the zeroth-order Bloch term to the same target DoFs. The higher-order mimetic stencils in \cite{paper2,paper3} replace the circulant coefficients while preserving the factorizations and commutativity used below. The matrix representation in 3D is assembled by Kronecker products between $\Dbf_s$ ($s=0,1$) and identity matrices:
\begin{align}
\label{display:dis:directional}
\Dcal_{s,1}
=I_{N^2}\otimes\Dbf_s,
\qquad
\Dcal_{s,2}=I_N\otimes\Dbf_s\otimes I_N,
\qquad
\Dcal_{s,3}=\Dbf_s\otimes I_{N^2}.
\end{align}
For a general primitive cell, the matrix $\Abf:=\begin{pmatrix}
    \vec a_1 & \vec a_2 & \vec a_3
\end{pmatrix}$ collects the lattice vectors, and $b_{ij}:=(\Abf^{-1})_{ij}$ denotes the $(i,j)$-th entry of its inverse. After mapping $\Omega$ to the unit cube, we define the transformed blocks as in \cite{jin2025kcpc,paper2}:
\begin{align}
\label{display:dis:general_directional}
\Dcal_i:=\sum_{j=1}^3b_{ji}\Dcal_{1,j}+\imagunit\wnkk_i\Dcal_{0,i},
\qquad i=1,2,3.
\end{align}
The discrete curl and divergence operators $\Acal$ and $\Bcal$ are defined by
\begin{align}
\label{display:dis:curl_div}
\Acal
:=
\begin{pmatrix}
0&-\Dcal_3&\Dcal_2\\
\Dcal_3&0&-\Dcal_1\\
-\Dcal_2&\Dcal_1&0
\end{pmatrix},
\qquad
\Bcal
:=
\begin{pmatrix}
\Dcal_1&\Dcal_2&\Dcal_3
\end{pmatrix}.
\end{align}

\begin{prop}\label{prop:dis:ab_prop}
    The discrete matrices $\Acal$ and $\Bcal$ satisfy  
    \begin{align}
        \label{display:dis:ba0}        \Bcal\Acal=0,\quad \Acal\Acal^\conjT+\Bcal^\conjT\Bcal=\diag(\Lcal,\Lcal,\Lcal),
    \end{align}
    where
    \begin{align}
        \label{display:dis:laplacian}        \Lcal:=\Bcal\Bcal^\conjT=\sum\limits_{i=1}^3\Dcal_i^\conjT\Dcal_i.
    \end{align}
\end{prop}

\begin{proof}
    The result follows from \cite[Proposition 2.2, Lemma 2.4]{paper2}.
\end{proof}

We define the diagonal matrices $\Mcal_0,\Mcal_1\in\Rbb^{3N^3\times3N^3}$ by
\begin{align}
\label{display:dis:m0_m1}
(\Mcal_1)_{\beta\beta}
:=
\begin{cases}
1,&\text{the $\beta$-th edge DoF lies in $\Omega_1$},\\
0,&\text{otherwise},
\end{cases}\qquad
\Mcal_0:=I_{3N^3}-\Mcal_1,\quad\beta=1,\ldots,3N^3.
\end{align}
The matrices $\Mcal_0$ and $\Mcal_1$ are complementary orthogonal projectors. The discrete inverse dielectric matrix is 
\begin{align}\label{display:dis:mcal}
    \Mcal_\eta(\omega):=\Mcal_0+\eta(\omega)\Mcal_1.
\end{align}
The discretization of \eqref{display:dis:model} gives the constrained nonlinear matrix eigenproblem
\begin{align}
\label{display:dis:form_nep_matrix}
\left\{
\begin{aligned}
\left(\Acal\Mcal_\eta(\omega)\Acal^\conjT\right)\vec H_h
&=\omega^2\vec H_h,\\
\Bcal\vec H_h&=0.
\end{aligned}
\right.
\end{align}
\subsection{Matrix-Free Formulation}
\label{sse:disc:mf}

The unitary one-dimensional DFT matrix $\Fbf$ and its three-dimensional and vector-valued forms are defined by
\begin{align}
\label{display:dis:dft}
\begin{aligned}
(\Fbf)_{pq}
&=N^{-\frac12}\exp\left(\frac{2\pi\imagunit(p-1)(q-1)}{N}\right),\\
\Fcal
&=\Fbf\otimes\Fbf\otimes\Fbf,
\qquad
\Fcal_3=\diag(\Fcal,\Fcal,\Fcal).
\end{aligned}
\end{align}
The matrices in \eqref{display:dis:circulant} are circulant, so each directional block admits the diagonal factorization
\begin{align}
\label{display:dis:directional_factor}
\Dcal_i=\Fcal\Kcal_i\Fcal^\conjT,
\qquad i=1,2,3.
\end{align}
Here $\Kcal_i$ is diagonal. The discrete matrices $\Acal$ and $\Bcal$ satisfy the factorizations
\begin{align}
\label{display:dis:ab_fac}
\Acal
=\Fcal_3\Kcal_\Acal\Fcal_3^\conjT,
\quad
\Bcal^\conjT\Bcal
=\Fcal_3\Kcal_\Bcal\Fcal_3^\conjT.
\end{align}
Here the Fourier symbols are
\begin{align}
\label{display:dis:symbols}
\Kcal_\Acal
:=
\begin{pmatrix}
0&-\Kcal_3&\Kcal_2\\
\Kcal_3&0&-\Kcal_1\\
-\Kcal_2&\Kcal_1&0
\end{pmatrix},
\quad
\Kcal_\Bcal
:=
\begin{pmatrix}
\Kcal_1^\conjT\\
\Kcal_2^\conjT\\
\Kcal_3^\conjT
\end{pmatrix}
\begin{pmatrix}
\Kcal_1&\Kcal_2&\Kcal_3
\end{pmatrix}.
\end{align}
In particular, $\Bcal\Acal=0$ of \eqref{display:dis:ba0} reduces to
\begin{align}
\label{display:dis:matrixfree_ba0}
\Kcal_\Bcal\Kcal_\Acal=0.
\end{align}
With $\vec x=\Fcal_3^\conjT\vec H_h$, the constrained problem becomes
\begin{align}
\label{display:dis:form_nep_matrixfree}
\left\{\begin{aligned}
\Kcal_0(\omega)\vec x&=\omega^2\vec x,
\\
\Kcal_\Bcal\vec x&=0,
\end{aligned}
\right.\quad \text{where }\Kcal_0(\omega)
:=
\Kcal_\Acal\Fcal_3^\conjT
\Mcal_\eta(\omega)
\Fcal_3\Kcal_\Acal^\conjT.
\end{align}
In practice, we store the nonzero elements of $\Kcal_\Acal$, $\Kcal_\Bcal$, and corresponding integer indices of $\Mcal_1$. Each application of $\Kcal_0(\omega)$ uses elementwise multiplication and one forward and one backward 3D FFT. For a block of $m$ vectors $\vec X\in\Cbb^{3N^3\times m}$, the dominant runtime cost is $O(mN^3\log N)$.

We conclude this section with the derivatives of $\Kcal_0(\omega)$. The material derivatives satisfy
\begin{align}
\label{display:dis:material_derivative}
\Mcal_\eta^{(\ell)}(\omega):=\frac{d^\ell\Mcal_\eta(\omega)}{d\omega^\ell}
=\eta^{(\ell)}(\omega)\Mcal_1,
\qquad \ell=1,2,\ldots.
\end{align}
The superscript ``$\ell$'' denotes the derivative order. Consequently, the derivatives of $\Kcal_0(\omega)$ are
\begin{align}
\label{display:dis:maxwell_derivatives}
\Kcal_0^{(\ell)}(\omega)
=\eta^{(\ell)}(\omega)\Kcal_\Acal\Fcal_3^\conjT
\Mcal_1
\Fcal_3\Kcal_\Acal^\conjT,
\qquad \ell=1,2,\ldots.
\end{align}
In a matrix-free implementation, $\Mcal_\eta(\omega)$ and its derivatives are applied by elementwise multiplications on the edge indices. We denote the frequency-independent term by
\begin{align}
\label{display:dis:q_block}
\Qcal
:=
\Kcal_\Acal\Fcal_3^\conjT
\Mcal_1
\Fcal_3\Kcal_\Acal^\conjT.
\end{align}
The derivatives therefore satisfy
$\Kcal_0^{(\ell)}(\omega)=\eta^{(\ell)}(\omega)\Qcal$ for $\ell=1,2,\ldots$.

\section{Kernel Compensation Formulation}
\label{se:kc}

The kernel compensation method adds a penalty term with coefficient $\gamma\geq0$ to $\Kcal_0(\omega)$:
\begin{align}
\label{display:kc:form}
\Kcal_\gamma(\omega)\vec x=\omega^2\vec x,\qquad
\Kcal_\gamma(\omega):=\Kcal_0(\omega)+\gamma\Kcal_\Bcal.
\end{align}
The value $\gamma=0$ recovers the original matrix $\Kcal_0(\omega)$.
Since the compensation term is independent of $\omega$, the derivatives satisfy
$$
\Kcal_\gamma^{(\ell)}(\omega)=\Kcal_0^{(\ell)}(\omega)=\eta^{(\ell)}(\omega)\Qcal,\qquad\ell=1,2,\ldots.
$$

In this section, we rebuild the kernel compensation formulation for the non-Hermitian Maxwell eigenproblem through unitary diagonalization, then extend the theory to nonlinear cases. In Subsection~\ref{sse:kc:preconditioning}, we prove that the preconditioning in \cite{paper2,paper3} remains efficient for nonlinear, non-Hermitian cases.

\subsection{Space Decomposition}
\label{sse:kc:spaces}

The following lemma gives the kernel of $\Kcal_0(\omega)$.
\begin{lem}
\label{lem:kc:ker_ama}
The matrices $\Acal$ and $\Mcal$ have compatible dimensions. If $\Mcal+\Mcal^\conjT$ is Hermitian positive definite (HPD), then the following equality holds:
$$
\ker(\Acal^\conjT\Mcal\Acal)=\ker(\Acal).
$$
\end{lem}
\begin{proof}
The inclusion $\ker(\Acal)\subset\ker(\Acal^\conjT\Mcal\Acal)$ follows directly. If $\Acal^\conjT\Mcal\Acal\vec x=\vec0$, then
$$
0=2\re\left(\vec x^\conjT\Acal^\conjT\Mcal\Acal\vec x\right)
=\left(\Acal\vec x\right)^\conjT\left(\Mcal+\Mcal^\conjT\right)\left(\Acal\vec x\right).
$$
The positive definiteness of $\Mcal+\Mcal^\conjT$ implies $\Acal\vec x=\vec0$, which proves the reverse inclusion.
\end{proof}

We introduce the discrete curl image space $\Vcal_\Acal$, the discrete divergence image space $\Vcal_\Bcal$, and the discrete harmonic space $\Vcal_0$ by
\begin{align}
\label{display:kc:spaces}
\Vcal_\Acal
:=\image(\Kcal_\Acal),
\qquad
\Vcal_\Bcal
:=\image(\Kcal_\Bcal),
\qquad
\Vcal_0
:=\ker(\Kcal_\Acal^\conjT)\cap\ker(\Kcal_\Bcal).
\end{align}
The following proposition gives the discrete space decomposition and identifies the kernel and image of $\Kcal_\gamma(\omega)$.
\begin{prop}
\label{prop:kc:kernel_image}
The spaces in \eqref{display:kc:spaces} are mutually orthogonal and satisfy
\begin{align}
\label{display:kc:c3n3_decomp}
\Cbb^{3N^3}=\Vcal_\Acal\oplus\Vcal_\Bcal\oplus\Vcal_0.
\end{align}
If $\re\eta(\omega)>0$, the kernels and images satisfy
\begin{align}
\label{display:kc:ker_im_pnt}
\begin{aligned}
\ker(\Kcal_\gamma(\omega))
=\ker(\Kcal_\gamma(\omega)^\conjT)
&=
\begin{cases}
\ker(\Kcal_\Acal^\conjT),&\gamma=0,\\
\Vcal_0,&\gamma>0.
\end{cases}\\
\image(\Kcal_\gamma(\omega))
=\image(\Kcal_\gamma(\omega)^\conjT)
&=
\begin{cases}
\Vcal_\Acal,&\gamma=0,\\
\Vcal_\Acal\oplus\Vcal_\Bcal,&\gamma>0.
\end{cases}
\end{aligned}
\end{align}
\end{prop}
\begin{proof}
Due to the range--kernel orthogonality, the following equalities hold:
$$
\Vcal_\Acal^\perp=\ker(\Kcal_\Acal^\conjT),
\qquad
\Vcal_\Bcal^\perp=\ker(\Kcal_\Bcal).
$$
The identities in \eqref{display:dis:matrixfree_ba0} imply $\Vcal_\Acal\perp\Vcal_\Bcal$. The harmonic space therefore satisfies
$$
\Vcal_0
=\Vcal_\Acal^\perp\cap\Vcal_\Bcal^\perp
=(\Vcal_\Acal\oplus\Vcal_\Bcal)^\perp,
$$
which proves the orthogonal decomposition in \eqref{display:kc:c3n3_decomp}.

The Hermitian part of $\Mcal_\eta(\omega)$ satisfies
$$
\frac{\Mcal_\eta(\omega)+\Mcal_\eta(\omega)^\conjT}{2}
=\Mcal_0+\re\eta(\omega)\Mcal_1
\succ0
$$
if $\re\eta(\omega)>0$ according to \eqref{display:dis:m0_m1} and \eqref{display:dis:mcal}. Lemma~\ref{lem:kc:ker_ama} then implies
$$
\ker(\Kcal_0(\omega))
=\ker(\Kcal_0(\omega)^\conjT)
=\ker(\Kcal_\Acal^\conjT).
$$
For $\gamma>0$ and every $\vec x\in\Cbb^{3N^3}$, the real part of the
quadratic form satisfies
$$
\re\left(\vec x^\conjT\Kcal_\gamma(\omega)\vec x\right)
=
\left(\Fcal_3\Kcal_\Acal^\conjT\vec x\right)^\conjT
\left(\Mcal_0+\re\eta(\omega)\Mcal_1\right)
\left(\Fcal_3\Kcal_\Acal^\conjT\vec x\right)
+\gamma\vec x^\conjT\Kcal_\Bcal\vec x.
$$
Both terms on the right-hand side are nonnegative. The conditions for
their sum to vanish yield
$$
\Kcal_\gamma(\omega)\vec x=\vec0 \quad\Leftrightarrow\quad
\Kcal_\Acal^\conjT\vec x=\vec0,\quad
\Kcal_\Bcal\vec x=\vec0,
$$
leading to $\ker(\Kcal_\gamma(\omega))=\Vcal_0$ for $\gamma>0$. The same argument gives $\ker(\Kcal_\gamma(\omega)^\conjT)=\Vcal_0$ because the adjoint has exactly the same Hermitian part.
Finally, $\image(\Kcal_\gamma(\omega))=\ker(\Kcal_\gamma(\omega)^\conjT)^\perp$ and $\image(\Kcal_\gamma(\omega)^\conjT)=\ker(\Kcal_\gamma(\omega))^\perp$ give the image formulas in \eqref{display:kc:ker_im_pnt}.
\end{proof}

The kernel of the discrete Laplacian determines the dimensions of these spaces.
For the simple cubic lattice $\Abf=I_3$, we assume $\max_{i=1,2,3}|\wnkk_i|\leq\pi$ and $h<\frac{1}{\pi}$.
\begin{prop}
\label{prop:kc:dim_lv}
Under these assumptions, the discrete Laplacian $\Lcal$ satisfies
\begin{align}
\label{display:kc:dim_l}
\dim\ker(\Lcal)=
\begin{cases}
0,&\wnk\ne\vec0,\\
1,&\wnk=\vec0,
\end{cases}
\end{align}
and the dimensions of $\Vcal_\Acal$, $\Vcal_\Bcal$, and $\Vcal_0$ are
\begin{align}
\label{display:kc:dim_v}
(\dim\Vcal_\Acal,\dim\Vcal_\Bcal,\dim\Vcal_0)
=
\begin{cases}
(2N^3,N^3,0),&\wnk\ne\vec0,\\
(2N^3-2,N^3-1,3),&\wnk=\vec0.
\end{cases}
\end{align}
\end{prop}
\begin{proof}
The dimension of $\ker(\Lcal)$ in \eqref{display:kc:dim_l} follows from \cite[Lemma 2.4]{paper2} for $\wnk\ne\vec0$; at $\wnk=\vec0$, $\ker(\Lcal)$ is spanned by the constant Fourier mode. The identities in \eqref{display:dis:ba0} and \eqref{display:dis:ab_fac} give
$$
\Kcal_\Acal\Kcal_\Acal^\conjT+\Kcal_\Bcal
=\Fcal_3^\conjT\diag(\Lcal,\Lcal,\Lcal)\Fcal_3.
$$
Therefore, $\Vcal_0=\ker(\Kcal_\Acal\Kcal_\Acal^\conjT+\Kcal_\Bcal)$ implies $\dim\Vcal_0=3\dim\ker(\Lcal)$. Since $\Bcal\Bcal^\conjT=\Lcal$, the rank satisfies
\begin{align}
\label{display:kc:rank_count}
\dim\Vcal_\Bcal=\rank(\Kcal_\Bcal)
=\rank(\Bcal^\conjT\Bcal)
=\rank(\Lcal)
=N^3-\dim\ker(\Lcal).
\end{align}
The orthogonal decomposition \eqref{display:kc:c3n3_decomp} then gives $\dim\Vcal_\Acal=3N^3-\dim\Vcal_\Bcal-\dim\Vcal_0$, which proves \eqref{display:kc:dim_v}.
\end{proof}
The same dimension formulas hold for other lattices whenever \eqref{display:kc:dim_l} holds. The Fourier symbols satisfy this condition for the face centered cubic and body centered cubic grids and wave vectors used in Section~\ref{se:num}.

\subsection{Unitary Block Diagonalization}
\label{sse:kc:unitary}

The following lemma identifies the spectrum retained by compression to the image space.
\begin{lem}
\label{lem:kc:image_compression}
The matrix $\Acal$ is square, and the columns of $\Ucal$ form an orthonormal basis of $\image(\Acal)$. The nonzero spectra agree, including algebraic multiplicities:
\begin{align}
\label{display:kc:image_spectrum}
\sigma(\Ucal^\conjT\Acal\Ucal)\setminus\{0\}
=\sigma(\Acal)\setminus\{0\}.
\end{align}
The matrix $\Ucal^\conjT\Acal\Ucal$ is nonsingular if and only if
$\image(\Acal)\cap\ker(\Acal)=\{\vec0\}$.
\end{lem}
\begin{proof}
Since $\Ucal\Ucal^\conjT\Acal=\Acal$, a unitary completion of $\Ucal$ transforms $\Acal$ into an upper triangular block matrix with diagonal blocks $\Ucal^\conjT\Acal\Ucal$ and $0$.
This proves \eqref{display:kc:image_spectrum} and preserves the nonzero algebraic multiplicities.
The identity $\Acal\Ucal=\Ucal(\Ucal^\conjT\Acal\Ucal)$ also gives
$\Ucal\ker(\Ucal^\conjT\Acal\Ucal)=\image(\Acal)\cap\ker(\Acal)$,
which proves the last assertion.
\end{proof}

We select arbitrary orthonormal bases of $\Vcal_\Acal$, $\Vcal_\Bcal$, and $\Vcal_0$, and denote the corresponding basis matrices by $\Ucal_\Acal$, $\Ucal_\Bcal$, and $\Ucal_0$, respectively. Each matrix has $3N^3$ rows and as many columns as the dimension of its space. Proposition~\ref{prop:kc:kernel_image} implies that the matrix
\begin{align}
\label{display:kc:unitary_bases}
\Ucal:=\begin{pmatrix}\Ucal_\Acal&\Ucal_\Bcal&\Ucal_0\end{pmatrix}
\end{align}
is unitary. Moreover, $\Ucal_0$ is omitted when $\Vcal_0=\{\vec0\}$. The restrictions of $\Kcal_0(\omega)$ and $\Kcal_\Bcal$ to their image spaces are defined by
\begin{align}
\label{display:kc:reduced_mat}
\Rcal_\Acal(\omega):=\Ucal_\Acal^\conjT\Kcal_0(\omega)\Ucal_\Acal,\quad
\Rcal_\Bcal:=\Ucal_\Bcal^\conjT\Kcal_\Bcal\Ucal_\Bcal.
\end{align}
For $\omega\in\Omega_{\mathrm{sol}}$, Lemma~\ref{lem:kc:image_compression} and Proposition~\ref{prop:kc:kernel_image} imply that $\Rcal_\Acal(\omega)$ and $\Rcal_\Bcal$ are nonsingular and retain all nonzero eigenvalues and algebraic multiplicities of $\Kcal_0(\omega)$ and $\Kcal_\Bcal$, respectively.
\begin{prop}
\label{prop:kc:unitary_decom_ab}
If $\omega\in\Omega_{\mathrm{sol}}$, then $\Rcal_\Acal(\omega)+\Rcal_\Acal(\omega)^\conjT$ and $\Rcal_\Bcal$ are HPD. The following factorizations hold:
\begin{align}
\label{display:kc:support}
\Kcal_0(\omega)=\Ucal_\Acal\Rcal_\Acal(\omega)\Ucal_\Acal^\conjT,\quad
\Kcal_\Bcal=\Ucal_\Bcal\Rcal_\Bcal\Ucal_\Bcal^\conjT.
\end{align}
\end{prop}
\begin{proof}
For any $\vec y\ne\vec0$, the vector $\Kcal_\Acal^\conjT\Ucal_\Acal\vec y$
is nonzero because $\Vcal_\Acal\perp\ker(\Kcal_\Acal^\conjT)$.
The associated quadratic form satisfies
$$
\re\left(\vec y^\conjT\Rcal_\Acal(\omega)\vec y\right)
=\left(\Fcal_3\Kcal_\Acal^\conjT\Ucal_\Acal\vec y\right)^\conjT
\left(\Mcal_0+\re\eta(\omega)\Mcal_1\right)
\left(\Fcal_3\Kcal_\Acal^\conjT\Ucal_\Acal\vec y\right)>0.
$$
Thus $\Rcal_\Acal(\omega)+\Rcal_\Acal(\omega)^\conjT$ is HPD. Since $\Ucal_\Acal\Ucal_\Acal^\conjT$ is the orthogonal projection onto $\Vcal_\Acal$, while Proposition~\ref{prop:kc:kernel_image} gives $\image(\Kcal_0(\omega))=\image(\Kcal_0(\omega)^\conjT)=\Vcal_\Acal$, we have
$$
\Ucal_\Acal\Ucal_\Acal^\conjT\Kcal_0(\omega)
=\Kcal_0(\omega)
=\Kcal_0(\omega)\Ucal_\Acal\Ucal_\Acal^\conjT.
$$
The two identities give
$$
\Kcal_0(\omega)=\Ucal_\Acal\Ucal_\Acal^\conjT\Kcal_0(\omega)\Ucal_\Acal\Ucal_\Acal^\conjT=\Ucal_\Acal\Rcal_\Acal(\omega)\Ucal_\Acal^\conjT.
$$
The same argument proves that $\Rcal_\Bcal$ is HPD and $\Kcal_\Bcal=\Ucal_\Bcal\Rcal_\Bcal\Ucal_\Bcal^\conjT$. 
\end{proof}

The factorizations \eqref{display:kc:support} imply the block diagonalization of the compensated matrix. We summarize the conclusion as the theorem below.
\begin{thm}
\label{thm:kc:block}
For $\omega\in\Omega_{\mathrm{sol}}$ and $\gamma\geq0$, the unitary matrix in \eqref{display:kc:unitary_bases} satisfies
\begin{align}
\label{display:kc:block_diag}
\Ucal^\conjT\Kcal_\gamma(\omega)\Ucal
=\begin{pmatrix}
\Rcal_\Acal(\omega)& & \\
 &\gamma\Rcal_\Bcal& \\
 & & 0
\end{pmatrix},
\end{align}
where the zero block at (3,3) is omitted when $\Vcal_0=\{\vec0\}$. 
\end{thm}
\begin{proof}
The result follows by substituting \eqref{display:kc:support} into \eqref{display:kc:form}.
\end{proof}

    For the Hermitian case, we prove in \cite{paper3,paper2} that the spectrum of $\Kcal_\gamma(\cdot)$ can be separated through the Courant-Fischer min-max principle, which does not hold for non-Hermitian matrices. Theorem~\ref{thm:kc:block} reveals that, when $\omega\in\Omega_{\mathrm{sol}}$ is fixed, the nonzero eigenvalues of $\Kcal_\gamma(\omega)$ consist of the nonzero eigenvalues of both $\Kcal_0(\omega)$ and $\gamma\Kcal_\Bcal$, which is compatible with the Hermitian case.

For the nonlinear, non-Hermitian case, the unitary block argument also proves the equivalence between the original problem and its restriction. 
\begin{prop}
\label{prop:kc:nep_decom_a}
For every $\omega\in\Omega_{\mathrm{sol}}$, the kernels satisfy
\begin{align}
\label{display:kc:nep_reduced}
\ker\left(\Kcal_0(\omega)-\omega^2I_{3N^3}\right)
=\Ucal_\Acal\ker\left(\Rcal_\Acal(\omega)-\omega^2I_{\dim\Vcal_\Acal}\right).
\end{align}
The correspondence $\vec x=\Ucal_\Acal\vec y$ has inverse $\vec y=\Ucal_\Acal^\conjT\vec x$, and every vector on the left satisfies $\Kcal_\Bcal\vec x=\vec0$.
\end{prop}
\begin{proof}
If $\Kcal_0(\omega)\vec x=\omega^2\vec x$, then $\vec x\in\image(\Kcal_0(\omega))=\Vcal_\Acal$. Due to $\Ucal_\Acal\Ucal_\Acal^\conjT\vec x=\vec x$, left multiplication by $\Ucal_\Acal^\conjT$ gives
$$
\Rcal_\Acal(\omega)(\Ucal_\Acal^\conjT\vec x)=\omega^2(\Ucal_\Acal^\conjT\vec x),
$$
thus $\vec y=\Ucal_\Acal^\conjT\vec x$ solves $\Rcal_\Acal(\omega)\vec y=\omega^2\vec y$. Conversely, \eqref{display:kc:support} shows that $\vec x=\Ucal_\Acal\vec y$ solves $\Kcal_0(\omega)\vec x=\omega^2\vec x$ for any $\vec y$ that solves $\Rcal_\Acal(\omega)\vec y=\omega^2\vec y$. 
\end{proof}

\begin{thm}
\label{thm:kc:nep}
If $\gamma>0$ satisfies
\begin{align}
\label{display:kc:gamma_cond_inthm}
\sup\left\{|\omega|^2:\omega\in\Omega_{\mathrm{sol}}\right\}
<\gamma\min\left\{\lambda>0:\lambda\in\sigma(\Kcal_\Bcal)\right\},
\end{align}
then the compensated problem \eqref{display:kc:form} and the original constrained problem \eqref{display:dis:form_nep_matrixfree} have the same eigenpairs for $\omega\in\Omega_{\mathrm{sol}}$.
\end{thm}
\begin{proof}
The block diagonalization \eqref{display:kc:block_diag} gives
\begin{align*}
    \det\left(\Kcal_\gamma(\omega)-\omega^2I_{3N^3}\right)
&=\det\left(\Ucal^\conjT\left(\Kcal_\gamma(\omega)-\omega^2I_{3N^3}\right)\Ucal\right)\\
    &=\det\left(\Rcal_\Acal(\omega)-\omega^2I_{\dim\Vcal_\Acal}\right)
  \det\left(\gamma\Rcal_\Bcal-\omega^2I_{\dim\Vcal_\Bcal}\right)
  (-\omega^2)^{\dim\Vcal_0}.
\end{align*}
The condition \eqref{display:kc:gamma_cond_inthm} makes $\gamma\Rcal_\Bcal-\omega^2I_{\dim\Vcal_\Bcal}$ nonsingular. Since $0\notin\Omega_{\mathrm{sol}}$, we have
$$
\det\left(\Kcal_\gamma(\omega)-\omega^2I_{3N^3}\right)=0
\quad\Leftrightarrow\quad
\det\left(\Rcal_\Acal(\omega)-\omega^2I_{\dim\Vcal_\Acal}\right)=0.
$$
Therefore $\Kcal_\gamma(\omega)\vec x=\omega^2\vec x$ and
$\det\left(\Rcal_\Acal(\omega)-\omega^2I_{\dim\Vcal_\Acal}\right)=0$
have the same eigenfrequencies. The $\Vcal_\Bcal$ and $\Vcal_0$ blocks of
$\Kcal_\gamma(\omega)-\omega^2I$ are nonsingular by
\eqref{display:kc:block_diag}, so every compensated eigenvector lies in
$\Vcal_\Acal$. Proposition~\ref{prop:kc:nep_decom_a} then gives the assertion.
\end{proof}

\subsection{Preconditioning}
\label{sse:kc:preconditioning}

The condition \eqref{display:kc:gamma_cond_inthm} places the positive
compensation eigenvalues above all $|\omega|^2$ in the search domain,
so the compensation block cannot contribute an eigenfrequency there.
In this paper, we use the penalty coefficient in \cite[Remark~2.6]{paper2}:
\begin{align}
\label{display:kc:pnt_coef}
\gamma
=
\begin{cases}
4\pi^2,&\wnk=\vec0\ \text{or}\ \normtwo{\wnk}>1,\\
4\pi^2\normtwo{\wnk}^{-2},&0<\normtwo{\wnk}\leq1.
\end{cases}
\end{align}
The Fourier symbols determine the smallest positive eigenvalue in
\eqref{display:kc:gamma_cond_inthm} for each grid and Bloch wave vector.
The chosen search domain must also satisfy this condition.

The frequency dependence of $\Kcal_\gamma(\omega)$ occurs only through $\Mcal_\eta(\omega)$. We follow \cite{paper2,paper3} and replace $\Mcal_\eta(\omega)$ by the identity to obtain the matrix $\Kcal_\Acal\Kcal_\Acal^\conjT+\gamma\Kcal_\Bcal$. The resulting matrix and $\Kcal_\gamma(\omega)$ have the common kernel $\Vcal_0$. We use the shifted preconditioner
\begin{align}
\label{display:kc:preconditioned_block}
\Kcal_\Pcal
:=
\Kcal_\Acal\Kcal_\Acal^\conjT
+\gamma\Kcal_\Bcal
+\delta I_{3N^3}.
\end{align}
In the implementation, the frozen matrix $\Kcal_\gamma(\omega)+\delta I$ and the preconditioner
use the same nonnegative shift
\begin{align}
\label{display:kc:common_shift}
\delta=
\begin{cases}
\pi^{-1},&\wnk=\vec0,\\
\normtwo{\wnk},&0<\normtwo{\wnk}\leq1,\\
0,&\normtwo{\wnk}>1.
\end{cases}
\end{align}
Under \eqref{display:kc:dim_l} and for $\gamma>0$,
Proposition~\ref{prop:kc:kernel_image} and \eqref{display:kc:common_shift}
imply that $\Kcal_\Pcal$ is HPD. Its $3\times3$ blocks are diagonal,
so its inverse is applied elementwise.
At $\wnk=\vec0$, the common shift replaces the three-dimensional harmonic
zero block in \eqref{display:kc:block_diag} by $\delta I_3$.
The shifted frozen matrix has the same eigenvectors as $\Kcal_\gamma$,
with every eigenvalue increased by $\delta$. The implementation orders
the shifted Ritz values by modulus and subtracts $\delta$ before the
frequency update. This ordering can differ from the ordering of the
unshifted eigenvalues by modulus.
The following proposition gives the spectrum of the preconditioned system.
\begin{prop}
\label{prop:kc:reduced_preconditioner}
For every $\omega\in\Omega_{\mathrm{sol}}$, the spectrum of the symmetrically preconditioned matrix satisfies
\begin{align}
\label{display:kc:sym_precond_spec}
\sigma\left(
\Kcal_\Pcal^{-\frac12}
\left(\Kcal_\gamma(\omega)+\delta I_{3N^3}\right)
\Kcal_\Pcal^{-\frac12}
\right)\subset
\left\{
1+t\left(\eta(\omega)-1\right):
0\leq t\leq1
\right\}.
\end{align}
\end{prop}
\begin{proof}
The kernel compensation formulation with shift $\delta$ can be rewritten as 
$$
\Kcal_\gamma(\omega)+\delta I_{3N^3}=\Kcal_\Pcal+(\eta(\omega)-1)\Qcal,
$$
where $\Qcal$ is defined in \eqref{display:dis:q_block}. The definition
of $\Qcal$ gives
$$
\Kcal_\Acal\Kcal_\Acal^\conjT-\Qcal=\Kcal_\Acal\Fcal_3^\conjT
\Mcal_0
\Fcal_3\Kcal_\Acal^\conjT.
$$
Since $\Mcal_0$ is positive semidefinite, the matrices satisfy
\begin{align}
\label{display:kc:matrix_ord}
0
\preceq
\Qcal
\preceq
\Kcal_\Acal\Kcal_\Acal^\conjT
\preceq
\Kcal_\Pcal.
\end{align}
The two-sided preconditioned matrix therefore satisfies
\begin{align}
\label{display:kc:sym_precond}
\Kcal_\Pcal^{-\frac12}
\left(\Kcal_\gamma(\omega)+\delta I_{3N^3}\right)
\Kcal_\Pcal^{-\frac12}
=
I_{3N^3}
+\left(\eta(\omega)-1\right)
\Kcal_\Pcal^{-\frac12}\Qcal\Kcal_\Pcal^{-\frac12}.
\end{align}
The Hermitian matrix
$\Kcal_\Pcal^{-\frac12}\Qcal\Kcal_\Pcal^{-\frac12}$
has eigenvalues in $[0,1]$ by \eqref{display:kc:matrix_ord}. Hence the matrix in \eqref{display:kc:sym_precond} has eigenvalues of the form
$1+t(\eta(\omega)-1)$ for some $t\in[0,1]$.
\end{proof}

\begin{rmk}
    In practice, we apply the left preconditioning $\Kcal_\Pcal^{-1}\left(\Kcal_\gamma(\omega)+\delta I_{3N^3}\right)$. Indeed, the relation
    $$
\Kcal_\Pcal^{\frac12}
\left[
\Kcal_\Pcal^{-1}
\left(\Kcal_\gamma(\omega)+\delta I_{3N^3}\right)
\right]
\Kcal_\Pcal^{-\frac12}=
\Kcal_\Pcal^{-\frac12}
\left(\Kcal_\gamma(\omega)+\delta I_{3N^3}\right)
\Kcal_\Pcal^{-\frac12}
$$
implies that the spectra of both preconditioned systems coincide. 
\end{rmk}

\begin{rmk}
    Neither $\Kcal_\gamma(\omega)$ nor the left-preconditioned matrix is normal in general. However, the symmetrically preconditioned matrix is normal because of \eqref{display:kc:sym_precond} and the fact that $\Kcal_\Pcal^{-\frac12}\Qcal\Kcal_\Pcal^{-\frac12}$ is Hermitian.
\end{rmk}

For every eigenvalue $\lambda$ in \eqref{display:kc:sym_precond_spec}, the assumption $0<\re\eta(\omega)<1$ in \eqref{display:dis:soldomain} gives
\begin{align}
\label{display:kc:preconditioned_bounds}
\begin{aligned}
\re\eta(\omega)
&\leq
\re\lambda
\leq
1,\\
\re\eta(\omega)
&\leq
|\lambda|
\leq
\max\{1,|\eta(\omega)|\}.
\end{aligned}
\end{align}
Since the two-sided matrix is normal, its condition number satisfies
\begin{align}
\label{display:kc:preconditioned_condition}
\kappa_2\left(
\Kcal_\Pcal^{-\frac12}
\left(\Kcal_\gamma(\omega)+\delta I_{3N^3}\right)
\Kcal_\Pcal^{-\frac12}
\right)
\leq
\frac{\max\{1,|\eta(\omega)|\}}
{\re\eta(\omega)}.
\end{align}
The same bound holds for the condition number of the left-preconditioned operator in the $\Kcal_\Pcal$-weighted norm. When $|\eta(\omega)|<1$, the upper bound in \eqref{display:kc:preconditioned_condition} becomes $\frac{1}{\re\eta(\omega)}$. In the Hermitian case
$\eta(\omega)=\varepsilon^{-1}(\omega)\in(0,1)$,
the complex line segment in \eqref{display:kc:sym_precond_spec} reduces to
$[\varepsilon^{-1}(\omega),1]$, and
\eqref{display:kc:preconditioned_condition} reduces to
$\kappa_2\leq\varepsilon(\omega)$. This recovers the preconditioning result in \cite{paper3}.

In Section~\ref{se:nep}, we take $\delta=0$ and assume $\wnk\ne\vec0$,
$\gamma>0$, and that \eqref{display:kc:dim_l} holds.
Proposition~\ref{prop:kc:kernel_image} then gives $\Vcal_0=\{\vec0\}$
and implies that $\Kcal_\Pcal$ is HPD.

\section{Solution of the Nonlinear Eigenproblem}
\label{se:nep}

This section presents the iterative solution of the nonlinear eigenproblem.
Subsection~\ref{sse:nep:linearized} introduces the linearized eigenproblem,
the regularity of its eigenpairs, and the Generalized Davidson eigensolver.
Subsection~\ref{sse:nep:newton} presents the quasi-Newton iteration and
the local regularity estimates used in its analysis.
Subsection~\ref{sse:nep:convergence} analyzes local convergence with exact
inner eigenpairs.

\subsection{Linearized Eigenproblem}
\label{sse:nep:linearized}
For each prescribed frequency $\omega\in\Omega_{\mathrm{sol}}$, we freeze the frequency dependence of $\Kcal_\gamma(\omega)$ and solve a non-Hermitian eigenproblem. The linearized problem is
\begin{align}
\label{display:nep:linearized_with_shift}
\Kcal_\gamma(\omega)\vec u_i(\omega)
=
\mu_i(\omega)\vec u_i(\omega),\quad i=1,\ldots,\nev,
\end{align}
where $\mu_i(\omega)$ is the $i$-th eigenvalue in increasing order of
magnitude, and $\vec u_i(\omega)$ is the corresponding normalized right
eigenvector. The number of desired eigenpairs is denoted by $\nev$.
Throughout this section, we assume that all desired eigenvalues are
algebraically simple:
\begin{align}
\label{display:nep:simple_eigenvalues}
\mu_i(\omega)
\text{ is algebraically simple},
\qquad
\omega\in\Omega_{\mathrm{sol}},
\quad
i=1,\ldots,\nev.
\end{align}
We also assume that the desired eigenvalues remain separated from the
remaining spectrum and that their increasing order of magnitude does not
change on $\Omega_{\mathrm{sol}}$. Hence, $\mu_i(\omega)$ consistently
denotes the $i$-th desired eigenvalue. Algebraic simplicity makes the
corresponding eigenspace one-dimensional, while unit normalization fixes
the magnitude of $\vec u_i(\omega)$.

The following theorem establishes the local holomorphy of $\mu_i$
and the regularity of its normalized eigenvectors.
\begin{thm}
\label{thm:nep:regularity}
We assume that $\eta$ is holomorphic near $\omega_0\in\Omega_{\mathrm{sol}}$
and that $\mu_i(\omega_0)$ is algebraically simple.
The eigenvalue admits a locally holomorphic simple branch $\mu_i(\omega)$.
The corresponding right and left eigenvectors can be chosen real analytic
with $\normtwo{\vec u_i(\omega)}=1$ and
\begin{align}
\label{display:nep:left_right_eig}
\Kcal_\gamma(\omega)^\conjT\vec v_i(\omega)
=\overline{\mu_i(\omega)}\vec v_i(\omega),
\qquad
\vec v_i(\omega)^\conjT\vec u_i(\omega)=1.
\end{align}
\end{thm}

\begin{proof}
The matrix $\Kcal_\gamma$ is holomorphic wherever $\eta$ is holomorphic.
Theorem~1 of \cite{greenbaum2020} therefore gives the locally holomorphic
simple eigenvalue branch $\mu_i$.
Theorem~2 of \cite{greenbaum2020} gives a holomorphic right eigenvector
$\vec x_i$ and a holomorphic left eigenvector row $\vec y_i^\conjT$
with $\vec y_i^\conjT\vec x_i=1$.
We normalize them by
$$
\vec u_i=\frac{\vec x_i}{\normtwo{\vec x_i}},
\qquad \vec v_i=\normtwo{\vec x_i}\,\vec y_i.
$$
The positive function $\normtwo{\vec x_i}$ and its reciprocal are real analytic
on a sufficiently small neighborhood, so $\vec u_i$ and $\vec v_i$
are real analytic and satisfy \eqref{display:nep:left_right_eig}.
\end{proof}

The following proposition expresses the exact eigenvalue derivative
through the right and left eigenvectors.
\begin{prop}
\label{prop:nep:mu_der}
Under the hypotheses of Theorem~\ref{thm:nep:regularity},
the derivative of the simple eigenvalue branch satisfies
\begin{align}
\label{display:nep:mu_i}
\der{\mu_i}(\omega)
=
\vec v_i(\omega)^\conjT
\der{\Kcal_\gamma}(\omega)
\vec u_i(\omega)
=
\der{\eta}(\omega)
\vec v_i(\omega)^\conjT
\Qcal
\vec u_i(\omega).
\end{align}
\end{prop}

\begin{proof}
The derivative formula in Theorem~1 of \cite{greenbaum2020} and
$\der{\Kcal_\gamma}=\der\eta\Qcal$ give \eqref{display:nep:mu_i}.
\end{proof}

The right eigenpairs are computed by a Generalized Davidson method with thick restart and soft locking. The method maintains an orthonormal search basis $\Zcal$ and solves a local eigenproblem defined by
\begin{align}
\label{display:nep:gd_projection}
\left(\Zcal^\conjT\Kcal_\gamma(\omega)\Zcal\right)\vec y_j=\theta_j\vec y_j,
\end{align}
during each iteration, where $\theta_j$ and $\vec y_j$ are an eigenvalue and the corresponding right unit eigenvector of the local matrix, respectively. Since $\Kcal_\gamma(\omega)$ is generally non-Hermitian, the vectors $\vec y_j$ need not be mutually orthogonal. The approximate right eigenvector $\vec x_j$ and residual $\vec r_j$ are defined by
\begin{align}
\label{display:nep:gd_residual}
\vec x_j:=\Zcal\vec y_j,\quad
\vec r_j:=\Kcal_\gamma(\omega)\vec x_j-\theta_j\vec x_j.
\end{align}
During each iteration, the search space $\Zcal$ is expanded by adding new directions computed from the preconditioned residuals $\Kcal_\Pcal^{-1}\vec r_j$ to the columns. We preserve part of the columns of $\Zcal$ and discard the rest when the number of columns of $\Zcal$ reaches a given maximum. The complete procedure is summarized in Algorithm~\ref{alg:nep:gd}.

\noindent\begin{minipage}{\linewidth}
\begin{breakablealgorithm}
\caption{Generalized Davidson Method with Thick Restart and Soft Locking for Solving the Exterior Eigenpairs of the Non-Hermitian Matrix $\Kcal_\gamma(\omega)$ with a Fixed Frequency $\omega$.}
\label{alg:nep:gd}
\begin{algorithmic}[1]
\Require The matrix $\Kcal_\gamma(\omega)$, the preconditioner $\Kcal_\Pcal^{-1}$, an initial block $\Zcal_0$ of rank $m\geq\nev$, tolerance $\tau_{\mathrm{GD}}>0$, maximum dimension of the search space $m_{\max}\geq2m$.
\Ensure Approximate right eigenpairs $(\mu_i(\omega),\vec u_i(\omega))$ for $i=1,\ldots,\nev$.
\State Set $\Zcal\gets\operatorname{orth}(\Zcal_0)$ and $\Wcal\gets\Kcal_\gamma(\omega)\Zcal$.
\While{the desired eigenpairs have not converged}
    \State Compute $\Hcal\gets\Zcal^\conjT\Wcal$ and the local eigenproblem $\Hcal\vec y_j=\theta_j\vec y_j$.
    \State The algorithm arranges $\theta_j$ in increasing order of magnitude and retains the first $m$ pairs.
    \State Compute $\vec x_j\gets\Zcal\vec y_j$ and residuals $\vec r_j\gets\Wcal\vec y_j-\theta_j\vec x_j$ for $j=1,\ldots,m$.
    \If{$\max\limits_{1\leq j\leq\nev}\normtwo{\vec r_j}\leq\tau_{\mathrm{GD}}$}
        \State Return $(\theta_i,\vec x_i)$ as $(\mu_i(\omega),\vec u_i(\omega))$ for $i=1,\ldots,\nev$.
    \EndIf
    \State Find the unconverged indices $\Ical\gets\{j\in\{1,\ldots,m\}:\normtwo{\vec r_j}>\tau_{\mathrm{GD}}\}$ and set $\Rcal\gets(\vec r_j)_{j\in\Ical}$.
    \If{$\rank(\Zcal)+|\Ical|>m_{\max}$}
        \State The thick restart sets $\Zcal\gets\operatorname{orth}([\vec x_1,\ldots,\vec x_m])$ and $\Wcal\gets\Kcal_\gamma(\omega)\Zcal$.
    \EndIf
    \State Compute the new directions by preconditioning the unconverged residuals: $\Tcal\gets\Kcal_\Pcal^{-1}\Rcal$.
    \State Apply $\Tcal\gets(I_{3N^3}-\Zcal\Zcal^\conjT)\Tcal$, then orthogonalize $\Tcal\gets\operatorname{orth}(\Tcal)$.
    \State Expand the space $\Zcal\gets[\Zcal,\Tcal]$ and $\Wcal\gets[\Wcal,\Kcal_\gamma(\omega)\Tcal]$.
\EndWhile
\end{algorithmic}
\end{breakablealgorithm}

\end{minipage}\par\medskip

Within the Davidson iteration, $\operatorname{orth}(\cdot)$ uses Cholesky QR
only for the restarted basis and projected preconditioned residual blocks.
Although forming the Gram matrix squares the block condition number,
these factorizations are numerically stable in our computations and
more efficient than Householder QR.

\begin{rmk}
\label{rmk:nep:jacobi_davidson}
The projected Jacobi--Davidson correction associated with
$\Kcal_\Pcal$ replaces lines 13--14 of Algorithm~\ref{alg:nep:gd} by
$$
\left(I_{3N^3}-\Zcal\Zcal^\conjT\right)
\Kcal_\Pcal
\left(I_{3N^3}-\Zcal\Zcal^\conjT\right)\Tcal
=-\Rcal,\qquad
\Zcal^\conjT\Tcal=0.
$$
Its solution is
$$
\Tcal=-\Kcal_\Pcal^{-1}\Rcal+\Kcal_\Pcal^{-1}\Zcal
\left(\Zcal^\conjT\Kcal_\Pcal^{-1}\Zcal\right)^{-1}
\Zcal^\conjT\Kcal_\Pcal^{-1}\Rcal.
$$
This requires computing and storing $\Kcal_\Pcal^{-1}\Zcal$.
Algorithm~\ref{alg:nep:gd} avoids this additional cost by applying
$\Kcal_\Pcal^{-1}$ only to the unconverged residuals $\Rcal$.
\end{rmk}

\begin{rmk}
\label{rmk:nep:storage_zw}
We preallocate $\Zcal$ and $\Wcal=\Kcal_\gamma\Zcal$ with
$m_{\max}$ columns, requiring $O(m_{\max}N^3)$ storage.
Each expansion applies $\Kcal_\gamma$ only to at most $m$ new directions
at a cost of $O(mN^3\log N)$, while residual preconditioning costs
$O(mN^3)$.
Projection and orthogonalization cost at most $O(m_{\max}^2N^3)$,
and the dense projected eigensolve costs $O(m_{\max}^3)$.
A restart recomputes the operator action on the retained basis.
The total cost per outer step depends on the numbers of expansions
and restarts.
\end{rmk}

\subsection{Quasi-Newton Iteration}
\label{sse:nep:newton}
The $i$-th characteristic frequency of
$\Kcal_\gamma(\omega)\vec x=\omega^2\vec x$ is a root of
\begin{align}
\label{display:nep:mu_omega2}
\omega^2-\mu_i(\omega)=0.
\end{align}
In \cite{paper3}, we apply Newton iteration to solve the corresponding
nonlinear equation in the Hermitian case. The exact Newton update for
\eqref{display:nep:mu_omega2} is
\begin{align}
    \label{display:nep:exact_newton}
    \omega_i^{(k+1)}\gets\omega_i^{(k)}-\frac{(\omega_i^{(k)})^2-\mu_i(\omega_i^{(k)})}{2\omega_i^{(k)}-\der{\mu_i}(\omega_i^{(k)})},
\end{align}
where $\omega_i^{(k)}$ denotes the approximate $i$-th frequency after
$k$ iterations. Since the formula for $\der\mu_i$ in \eqref{display:nep:mu_i} involves both left and
right eigenvectors, we replace $\der\mu_i(\omega)$ by $\vec u_i(\omega)^\conjT\der{\Kcal_\gamma}(\omega)\vec u_i(\omega)$ as an approximate derivative. We denote the right Rayleigh quotient with respect to $\Qcal$ as 
\begin{align}\label{display:nep:qi_def}
    q_i(\omega):=\vec u_i(\omega)^\conjT\Qcal\vec u_i(\omega).
\end{align}
The quasi-Newton update is defined by
\begin{align}
\label{display:nep:qnewton}
\omega_i^{(k+1)}
\gets
\omega_i^{(k)}
-
\frac{
(\omega_i^{(k)})^2-\mu_i(\omega_i^{(k)})
}{
2\omega_i^{(k)}-(\der\eta q_i)(\omega_i^{(k)})
}.
\end{align}
In general, $(\der\eta q_i)(\omega)$ differs from $\der{\mu_i}(\omega)$ for the
non-normal, non-Hermitian matrix $\Kcal_\gamma(\omega)$.

At frequencies with $|\omega_i^{(k)}|>10^{-10}$, we accept the quasi-Newton
candidate only when it and its inverse dielectric coefficient are finite
and the numerical safeguard conditions are satisfied:
\begin{align}\label{display:nep:safeguard}
\left|2\omega_i^{(k)}-(\der\eta q_i)(\omega_i^{(k)})\right|>10^{-30},
\qquad \re\eta(\omega_i^{(k+1)})\geq0.
\end{align}
Otherwise we set $\omega_i^{(k+1)}\gets\sqrt{\mu_i(\omega_i^{(k)})}$.
The local convergence analysis concerns frequencies in
$\Omega_{\mathrm{sol}}$; the numerical safeguard allows a wider set of candidates.
The nonlinear residual of each step is defined by
\begin{align}
\label{display:nep:qnewton_residual}
\vec r_i^{(k+1)}
:=
\Kcal_\gamma(\omega_i^{(k+1)})
\vec u_i(\omega_i^{(k)})
-
(\omega_i^{(k+1)})^2
\vec u_i(\omega_i^{(k)}).
\end{align}
The iteration terminates when
$\normtwo{\vec r_i^{(k+1)}}\leq\tau_{\mathrm{NEP}}$.
Algorithm~\ref{alg:nep:newton} summarizes the procedure of the quasi-Newton method.
The returned vector is computed at the preceding outer frequency;
\eqref{display:nep:qnewton_residual} tests the approximate eigenpair at
the accepted new frequency.
The complete-band and tolerance experiments use the inner tolerance
\begin{align}
\label{display:nep:inner_tolerance}
\tau_{\mathrm{GD}}=0.9\tau_{\mathrm{NEP}}.
\end{align}
Subsection~\ref{sse:num:histories} specifies the tolerances used for the
convergence histories.

\noindent\begin{minipage}{\linewidth}
\begin{breakablealgorithm}
\caption{A Quasi-Newton Method for Solving the Nonlinear Eigenproblem.}
\label{alg:nep:newton}
\begin{algorithmic}[1]
\Require Matrix-valued function $\Kcal_\gamma(\omega)$, initial frequency $\omega_1^{(0)}$, number of required eigenpairs $\nev$, and tolerances $\tau_{\mathrm{NEP}}$ and $\tau_{\mathrm{GD}}$.
\Ensure Approximate frequencies $\omega_i$ and corresponding right eigenvectors $\vec x_i$\\
for $i=1,\ldots,\nev$.
\For{$i=1,\ldots,\nev$}
    \For{$k=0,\ldots$}
        \State Compute the normalized $i$-th right eigenpair
        $(\mu_i(\omega_i^{(k)}),\vec u_i(\omega_i^{(k)}))$
        by the Generalized Davidson method in Algorithm~\ref{alg:nep:gd} with $i$ target eigenpairs and tolerance
        $\tau_{\mathrm{GD}}$.
        \State Compute
        $q_i(\omega_i^{(k)})
        \gets
        \vec u_i(\omega_i^{(k)})^\conjT
        \Qcal
        \vec u_i(\omega_i^{(k)})$.
        \State Apply \eqref{display:nep:qnewton}, or use $\omega_i^{(k+1)}\gets\sqrt{\mu_i(\omega_i^{(k)})}$ if the acceptance checks above fail.
        \State Compute the nonlinear residual $\vec r_i^{(k+1)}$ by
        \eqref{display:nep:qnewton_residual}. Break the loop if $\normtwo{\vec r_i^{(k+1)}}\leq\tau_{\mathrm{NEP}}$.
    \EndFor
    \State Set $\omega_i\gets\omega_i^{(k+1)}$ and $\vec x_i\gets \vec u_i(\omega_i^{(k)})$.
    \If{$i<\nev$}
        \State Set $\omega_{i+1}^{(0)}\gets\sqrt{\mu_{i+1}(\omega_i^{(k)})}$ using the additional Ritz value retained by the inner solve.
    \EndIf
\EndFor
\end{algorithmic}
\end{breakablealgorithm}
\end{minipage}

The retained block size $m>i$ in Algorithm~\ref{alg:nep:gd} increases with the
current target $i$. The additional Ritz value used to initialize the next
frequency need not satisfy the inner stopping tolerance.
The fixed-point step uses the principal square root. For the next-frequency
initialization, its sign is chosen closest to the preceding frequency.

Theorem~\ref{thm:nep:regularity} implies that
$q_i:\Rbb^2\to\Rbb$ and $\der\eta q_i:\Rbb^2\to\Cbb$
are locally real analytic, with $\omega=x+\imagunit y$ identified with $(x,y)$.
We write
$$
\nabla q_i:=(\p_xq_i,\p_yq_i),
$$
and define its product with $e\in\Cbb$ as
$$
\langle\nabla q_i,e\rangle:=(\p_xq_i)\re e+(\p_yq_i)\impart e.
$$
Since $q_i$ is real-valued, it can be holomorphic only if it is locally constant. For a general frequency $\omega$, we omit the argument ``$(\omega)$'' hereafter.
The superscript ``$*$'' along with the subscript ``$i$'' denotes the evaluation at $\omega_i^*$, for example:
$$
q_i^*=q_i(\omega_i^*),\quad (\der\eta)_i^*=\der\eta(\omega_i^*),\quad(\p_xq_i)^*=\p_xq_i(\omega_i^*).
$$

We further assume that
\begin{align}
\label{display:nep:assump_redeta}
\re\der\eta<0,\quad\text{for }\omega\in\Omega_{\mathrm{sol}}.
\end{align}
For real frequencies with $\varepsilon>0$ and $\der\varepsilon>0$,
this follows from $\der\eta=-\frac{\der\varepsilon}{\varepsilon^2}<0$.
Continuity preserves the strict inequality under sufficiently small
complex perturbations.
For $\re\omega>0$, the update denominator satisfies
\begin{align}
\label{display:nep:denominator_bound}
\left|2\omega-\der\eta q_i\right|
\geq\re\left(2\omega-\der\eta q_i\right)
=2\re\omega-(\re\der\eta)q_i>0.
\end{align}
This bound is uniform on a sufficiently small closed neighborhood
of each frequency in $\Omega_{\mathrm{sol}}$ with $\re\omega>0$. The safeguard
\eqref{display:nep:safeguard} checks the denominator numerically.

\subsection{Analysis of Convergence}
\label{sse:nep:convergence}

In this subsection, we analyze the generally linear local convergence
of the quasi-Newton iteration and estimate its linear factor through
the gradient of $q_i$.
Throughout this subsection, we assume $\gamma>0$.

\subsubsection{Expansion of the Error Term}
\label{ssse:nep:err_exp}

We analyze the update \eqref{display:nep:qnewton} with exact eigenpairs.
We denote the frequency error after $k$ iterations by $e_i^{(k)}:=\omega_i^{(k)}-\omega_i^*$.
The following proposition gives the local error expansion, with
$\rho_i^*$ evaluated at $\omega_i^*$.

\begin{prop}
\label{prop:nep:local}
We assume that $\omega_i^*\in\Omega_{\mathrm{sol}}$ has positive real part and satisfies
\begin{align}
\label{display:nep:simple_root}
(\omega_i^*)^2-\mu_i^*=0,
\qquad
2\omega_i^*-(\der{\mu_i})^*\ne0.
\end{align}
For $\omega_i^{(k)}$ sufficiently close to $\omega_i^*$, the error
$e_i^{(k)}$ satisfies
\begin{align}
\label{display:nep:err_exp_complete}
e_i^{(k+1)}=\rho_i^*e_i^{(k)}+O\left(|e_i^{(k)}|^2\right).
\end{align}
Here, the linear error factor $\rho_i$ is defined by
\begin{align}
\label{display:nep:rho}
\rho_i:=\frac{\der\mu_i-\der\eta q_i}{2\omega-\der\eta q_i}=
\frac{\der\eta(\vec v_i-\vec u_i)^\conjT\Qcal\vec u_i}
{2\omega-\der\eta q_i}.
\end{align}
\end{prop}

\begin{proof}
The holomorphy of $\mu_i$ and the root condition
\eqref{display:nep:simple_root} give
$$
(\omega_i^*+e_i^{(k)})^2-\mu_i(\omega_i^*+e_i^{(k)})
=\left(2\omega_i^*-(\der{\mu_i})^*\right)e_i^{(k)}
+O\left(|e_i^{(k)}|^2\right).
$$
The real analyticity of $q_i$ and the denominator bound
\eqref{display:nep:denominator_bound} also give
$$
\frac{1}{2(\omega_i^*+e_i^{(k)})-(\der\eta q_i)(\omega_i^*+e_i^{(k)})}
=\frac{1}{2\omega_i^*-(\der\eta)_i^*q_i^*}
+O\left(|e_i^{(k)}|\right).
$$
Substitution into \eqref{display:nep:qnewton} yields
\begin{align}
\label{display:nep:err_exp_step1}
\begin{aligned}
e_i^{(k+1)}
&=\left(1-\frac{2\omega_i^*-(\der{\mu_i})^*}
{2\omega_i^*-(\der\eta)_i^*q_i^*}\right)e_i^{(k)}
+O\left(|e_i^{(k)}|^2\right)\\
&=\rho_i^*e_i^{(k)}+O\left(|e_i^{(k)}|^2\right).
\end{aligned}
\end{align}
The last equality uses \eqref{display:nep:rho}.
Thus \eqref{display:nep:err_exp_step1} proves
\eqref{display:nep:err_exp_complete}.
\end{proof}
For the exact Newton update \eqref{display:nep:exact_newton}, the same
Taylor expansion with $\der\mu_i$ in the denominator cancels the linear
term and gives local quadratic convergence at a simple root.
For \eqref{display:nep:qnewton}, the expansion gives local linear
convergence when $0<|\rho_i^*|<1$ and at least quadratic convergence
when $\rho_i^*=0$.

The following proposition expresses the linear factor $\rho_i$
in terms of the gradient of $q_i$.
\begin{prop}
\label{prop:nep:rho_grad}
Under the assumptions of Theorem~\ref{thm:nep:regularity},
the factor $\rho_i$ satisfies
\begin{align}
\label{display:nep:rho_rep}
\rho_i
=\frac{\impart\eta}{2\omega-\der\eta q_i}
\left(\p_yq_i+\imagunit\p_xq_i\right).
\end{align}
As a result, the modulus $|\rho_i|$ equals
\begin{align}
\label{display:nep:rho_grad_abs}
|\rho_i|
=\frac{|\impart\eta|}{\left|2\omega-\der\eta q_i\right|}
|\nabla q_i|.
\end{align}
\end{prop}
\begin{proof}
The real and imaginary parts of
$\mu_i=\vec u_i^\conjT\Kcal_\Pcal\vec u_i+(\eta-1)q_i$ give
\begin{align}
\label{display:nep:qi_prop}
q_i\leq\frac{\re\mu_i}{\re\eta},\quad
\impart\mu_i=q_i\impart\eta,
\end{align}
respectively. The Cauchy--Riemann equations give
$$
\der\mu_i
=\p_x(\re\mu_i)+\imagunit\p_x(\impart\mu_i)
=\p_y(\impart\mu_i)+\imagunit\p_x(\impart\mu_i).
$$
The relation $\impart\mu_i=q_i\impart\eta$ yields
\begin{align}
\label{display:nep:mu_q_der}
\der{\mu_i}
=\der\eta q_i
+\impart\eta\left(\p_yq_i+\imagunit\p_xq_i\right).
\end{align}
The formulas \eqref{display:nep:mu_q_der} and \eqref{display:nep:rho}
give \eqref{display:nep:rho_rep}.
\end{proof}

The representation \eqref{display:nep:rho_grad_abs} reduces the estimate
of $|\rho_i|$ to a bound on $|\nabla q_i|$.
The next subsection derives this bound from the relative spectral gap
of the Hermitian part.

\subsubsection{Gradient Bound from Spectral Regularity}
\label{ssse:nep:eigvec_diff}

For the real-valued function $q_i=q_i(x,y)$,
\eqref{display:nep:mu_i} and \eqref{display:nep:mu_q_der} give
\begin{align}
\label{display:nep:gradient_defect}
|\impart\eta|\,|\nabla q_i|
=|\der\eta|\,\left|(\vec v_i-\vec u_i)^\conjT\Qcal\vec u_i\right|.
\end{align}
We bound $\left|(\vec v_i-\vec u_i)^\conjT\Qcal\vec u_i\right|$ by comparing the
eigenpairs of $\Kcal_\gamma$ with those of its Hermitian part
\begin{align}
\label{display:nep:Hgamma}
\Hcal_\gamma:=\frac{\Kcal_\gamma+\Kcal_\gamma^\conjT}{2}.
\end{align}
Under the standing assumptions, \eqref{display:kc:matrix_ord} and
$0<\re\eta<1$ give
$$
\Hcal_\gamma=\Kcal_\Pcal-(1-\re\eta)\Qcal
\succeq(\re\eta)\Kcal_\Pcal\succ0.
$$
We use an orthonormal eigenbasis $\{\vec w_j\}$ with
$\Hcal_\gamma\vec w_j=\nu_j\vec w_j$ and the $\nu_j$ in nondecreasing order. The next proposition bounds the difference of $\mu_i$ and $\nu_i$ using the relative Bauer--Fike
estimate \cite{eisenstat1998relative}. For any HPD matrix
$\Hcal$ and perturbation matrix $\Ecal$, each eigenvalue $\lambda\in\sigma(\Hcal+\Ecal)$ satisfies
\begin{align}
\label{display:nep:relative_perturbation}
\min_{\nu\in\sigma(\Hcal)}\left|1-\frac{\lambda}{\nu}\right|
\leq\normtwo{\Hcal^{-\frac12}\Ecal\Hcal^{-\frac12}}.
\end{align}
Indeed, $I-\lambda\Hcal^{-1}+\Hcal^{-\frac12}\Ecal\Hcal^{-\frac12}$
is singular. The smallest singular value of the normal matrix
$I-\lambda\Hcal^{-1}$ is therefore bounded by the perturbation norm,
which gives \eqref{display:nep:relative_perturbation}.
\begin{prop}
\label{prop:nep:relative_localization}
We assume that $\nu_i$ is simple and that its relative gap satisfies
\begin{align}
\label{display:nep:assump_gap}
\Delta_i:=\min_{\nu_j\ne\nu_i}\frac{|\nu_i-\nu_j|}{\nu_i+\nu_j}
>\frac{|\impart\eta|}{\re\eta}.
\end{align}
Then the $i$-th eigenvalue $\mu_i$ of $\Kcal_\gamma$ in order of increasing
modulus is algebraically simple and satisfies
\begin{align}
\label{display:nep:relative_bound}
\min_j\left|1-\frac{\mu_i}{\nu_j}\right|
=\left|1-\frac{\mu_i}{\nu_i}\right|\leq\frac{|\impart\eta|}{\re\eta}.
\end{align}
\end{prop}
\begin{proof}
We set $\Ecal:=\Kcal_\gamma-\Hcal_\gamma=\imagunit(\impart\eta)\Qcal$
and consider
the homotopy $\Hcal_\gamma+t\Ecal$, $t\in[0,1]$. Since
$\Hcal_\gamma-(\re\eta)\Qcal=\Kcal_\Pcal-\Qcal\succeq0$ and $\Qcal\succeq0$,
\eqref{display:nep:relative_perturbation} gives, for each eigenvalue
$\lambda\in\sigma(\Hcal_\gamma+t\Ecal)$,
\begin{align}
\label{display:nep:relative_perturbation_path}
\min_j\left|1-\frac{\lambda}{\nu_j}\right|
\leq t|\impart\eta|\normtwo{\Hcal_\gamma^{-\frac12}\Qcal\Hcal_\gamma^{-\frac12}}
\leq\frac{t|\impart\eta|}{\re\eta}.
\end{align}
Since the $\nu_j$ are positive, $0<\Delta_i<1$.
The gap condition \eqref{display:nep:assump_gap} gives
$|\nu_i-\nu_j|\geq\Delta_i(\nu_i+\nu_j)$ for $j\ne i$. This yields
\begin{align}
\label{display:nep:gap_separation}
(1+\Delta_i)\nu_j\leq(1-\Delta_i)\nu_i\ \ (j<i),
\quad
(1-\Delta_i)\nu_j\geq(1+\Delta_i)\nu_i\ \ (j>i).
\end{align}
For $|\lambda|\in\left[(1-\Delta_i)\nu_i,(1+\Delta_i)\nu_i\right]$,
\eqref{display:nep:gap_separation} implies $|\lambda|/\nu_j\geq1+\Delta_i$
for $j<i$ and $|\lambda|/\nu_j\leq1-\Delta_i$ for $j>i$. Hence
\begin{align}
\label{display:nep:gap_relative_dist}
\left|1-\frac{\lambda}{\nu_j}\right|
\geq\left|1-\frac{|\lambda|}{\nu_j}\right|\geq\Delta_i,\quad j\ne i.
\end{align}
On either boundary circle, the same lower bound holds for $j=i$.
The strict inequality in \eqref{display:nep:assump_gap} therefore
excludes homotopy eigenvalues on either circle by
\eqref{display:nep:relative_perturbation_path}.

At $t=0$, the open disks bounded by the inner and outer circles contain
$i-1$ and $i$ eigenvalues, respectively, counted with algebraic multiplicity.
By continuity of the eigenvalues, the counts remain constant for all $t$.
The annulus therefore contains one algebraically simple eigenvalue
$\mu_i(t)$, which is the $i$-th eigenvalue in order of increasing modulus.
In particular,
\begin{align}
\label{display:nep:modulus_localization}
(1-\Delta_i)\nu_i<|\mu_i(t)|<(1+\Delta_i)\nu_i,\quad t\in[0,1].
\end{align}
By \eqref{display:nep:gap_relative_dist}, the minimizing index in
\eqref{display:nep:relative_perturbation_path} is $i$ along this branch.
At $t=1$, this gives \eqref{display:nep:relative_bound}.
\end{proof}

Since $\normtwo{\vec u_i}=1$ and $\vec v_i^\conjT\vec u_i=1$, we have
\begin{align}
\label{display:nep:left_right_identity}
\normtwo{\vec v_i-\vec u_i}^2=\normtwo{\vec v_i}^2-1,
\end{align}
so we bound $\normtwo{\vec v_i}$ through the components along $\vec w_i$.
\begin{prop}
\label{prop:nep:target_component}
If $\nu_i$ is simple and \eqref{display:nep:assump_gap} holds, then
\begin{align}
\label{display:nep:wix_estimate}
|\vec w_i^\conjT\vec x|^2>\frac{2+\sqrt3}{4},\quad \vec x\in\left\{\vec u_i,\frac{\vec v_i}{\normtwo{\vec v_i}}\right\}.
\end{align}
\end{prop}
\begin{proof}
We consider the defect
$(\vec x^\conjT\Hcal_\gamma\vec x)(\vec x^\conjT\Hcal_\gamma^{-1}\vec x)-1$.
For either unit eigenvector, $\vec x^\conjT\Kcal_\gamma\vec x=\mu_i$.
Its real part gives $\vec x^\conjT\Hcal_\gamma\vec x=\re\mu_i$.
If $\impart\eta\ne0$, its imaginary part gives
$\vec x^\conjT\Qcal\vec x=q_i$. If $\impart\eta=0$, the matrix is Hermitian
and the normalized left and right eigenvectors coincide by simplicity,
so the same identity holds. The relation $\Kcal_\gamma=\Hcal_\gamma+\imagunit(\impart\eta)\Qcal$ gives
$$
\Kcal_\gamma^\conjT\Hcal_\gamma^{-1}\Kcal_\gamma
=\Kcal_\gamma\Hcal_\gamma^{-1}\Kcal_\gamma^\conjT=\Hcal_\gamma+|\impart\eta|^2\Qcal\Hcal_\gamma^{-1}\Qcal.
$$
When $\vec x=\vec u_i$, $\Kcal_\gamma\vec x=\mu_i\vec x$ gives
$$
|\mu_i|^2\vec x^\conjT\Hcal_\gamma^{-1}\vec x
=(\Kcal_\gamma\vec x)^\conjT\Hcal_\gamma^{-1}(\Kcal_\gamma\vec x)
=\re\mu_i+|\impart\eta|^2\vec x^\conjT\Qcal\Hcal_\gamma^{-1}\Qcal\vec x.
$$
For the left eigenvector, $\Kcal_\gamma^\conjT\vec x=\overline{\mu_i}\vec x$
gives the same equality. We substitute $\vec x^\conjT\Hcal_\gamma\vec x=\re\mu_i$
and $q_i=\vec x^\conjT\Qcal\vec x$ to obtain the $q_i$-dependent defect bound
\begin{align*}
(\re\mu_i)\vec x^\conjT\Hcal_\gamma^{-1}\vec x-1
&=\frac{|\impart\eta|^2}{|\mu_i|^2}
\left(\re\mu_i\,\vec x^\conjT\Qcal\Hcal_\gamma^{-1}\Qcal\vec x-q_i^2\right)\notag\\
&\leq\frac{|\impart\eta|^2q_i\left(\re\mu_i-(\re\eta)q_i\right)}
{(\re\eta)|\mu_i|^2}\\
&\leq\frac{|\impart\eta|^2}{4(\re\eta)^2}\frac{|\re\mu_i|^2}{|\mu_i|^2},
\end{align*}
where the first inequality uses $\Qcal\Hcal_\gamma^{-1}\Qcal\preceq\Qcal/(\re\eta)$,
which follows from $\Hcal_\gamma\succeq(\re\eta)\Qcal$. The second uses the
arithmetic--geometric mean inequality
$q_i\left(\re\mu_i-(\re\eta)q_i\right)\leq(\re\mu_i)^2/(4\re\eta)$. The defect therefore also satisfies
\begin{align}
\label{display:nep:defect_upper}
(\re\mu_i)\vec x^\conjT\Hcal_\gamma^{-1}\vec x-1
\leq\frac{|\impart\eta|^2}{4(\re\eta)^2}<\frac{\Delta_i^2}{4}.
\end{align}
We expand $\vec x$ in the eigenbasis $\{\vec w_j\}$, pair $(j,k)$ with
$(k,j)$, and retain the terms involving $i$ to obtain
\begin{align*}
(\re\mu_i)\vec x^\conjT\Hcal_\gamma^{-1}\vec x-1
&=\sum_{j,k}\left(\frac{\nu_j}{\nu_k}-1\right)
|\vec w_j^\conjT\vec x|^2|\vec w_k^\conjT\vec x|^2
\\ &=\sum_{j<k}\frac{(\nu_j-\nu_k)^2}{\nu_j\nu_k}
|\vec w_j^\conjT\vec x|^2|\vec w_k^\conjT\vec x|^2 \notag\\
&\geq|\vec w_i^\conjT\vec x|^2
\sum_{j\ne i}\frac{(\nu_j-\nu_i)^2}{\nu_j\nu_i}|\vec w_j^\conjT\vec x|^2,
\end{align*}
where the last step keeps only the terms $j=i$ or $k=i$, all of which are
nonnegative. For $j\ne i$, the assumption \eqref{display:nep:assump_gap} gives
$$
\frac{(\nu_j-\nu_i)^2}{\nu_j\nu_i}
\geq\Delta_i^2\frac{(\nu_j+\nu_i)^2}{\nu_j\nu_i}
\geq4\Delta_i^2.
$$
This coefficient bound implies
\begin{align}
\label{display:nep:defect_lower}
(\re\mu_i)\vec x^\conjT\Hcal_\gamma^{-1}\vec x-1
\geq4\Delta_i^2|\vec w_i^\conjT\vec x|^2
\left(1-|\vec w_i^\conjT\vec x|^2\right).
\end{align}
Together, \eqref{display:nep:defect_upper} and \eqref{display:nep:defect_lower} give
\begin{align}
\label{display:nep:product_bound}
    |\vec w_i^\conjT\vec x|^2
\left(1-|\vec w_i^\conjT\vec x|^2\right)
<\frac{1}{16}.
\end{align}
The relative bound \eqref{display:nep:relative_bound} and the strict
gap condition \eqref{display:nep:assump_gap} give
$$
\left|\re\mu_i-\nu_i\right|\leq|\mu_i-\nu_i|
\leq\frac{\nu_i|\impart\eta|}{\re\eta}<\Delta_i\nu_i.
$$
The real-part localization and its spectral expansion are
\begin{align}
\label{display:nep:remu_loc}
(1-\Delta_i)\nu_i<\re\mu_i<(1+\Delta_i)\nu_i,\quad\re\mu_i=\sum_j\nu_j|\vec w_j^\conjT\vec x|^2.
\end{align}
The defect also admits the expansion
\begin{align*}
    (\re\mu_i)\vec x^\conjT\Hcal_\gamma^{-1}\vec x-1&=\sum\limits_j\left(\frac{\re\mu_i}{\nu_j}-1\right)|\vec w_j^\conjT\vec x|^2\\
    &=\sum\limits_j\left(\frac{\re\mu_i}{\nu_j}+\frac{\nu_j}{\re\mu_i}-2\right)|\vec w_j^\conjT\vec x|^2\\
    &=\sum\limits_j\frac{(\re\mu_i-\nu_j)^2}{(\re\mu_i)\nu_j}|\vec w_j^\conjT\vec x|^2,
\end{align*}
The real-part localization \eqref{display:nep:remu_loc} and the separation
\eqref{display:nep:gap_separation} imply that
$(\re\mu_i)/\nu_j>1+\Delta_i$ for $j<i$ and
$(\re\mu_i)/\nu_j<1-\Delta_i$ for $j>i$.
Thus $\left|1-\frac{\re\mu_i}{\nu_j}\right|\geq\Delta_i$ for $j\ne i$.
The corresponding coefficient satisfies
$$
\frac{(\nu_j-\re\mu_i)^2}{(\re\mu_i)\nu_j}
\geq
\frac{\left|1-\frac{\re\mu_i}{\nu_j}\right|^2}
{1+\left|1-\frac{\re\mu_i}{\nu_j}\right|}\geq \frac{\Delta_i^2}{1+\Delta_i}.
$$
The $j=i$ term in the spectral sum is nonnegative, and the weights for
$j\ne i$ sum to $1-|\vec w_i^\conjT\vec x|^2$. The resulting lower
bound, together with \eqref{display:nep:defect_upper} and
\eqref{display:nep:assump_gap}, gives
$$
\frac{\Delta_i^2}{1+\Delta_i}\left(1-|\vec w_i^\conjT\vec x|^2\right)
\leq
(\re\mu_i)\vec x^\conjT\Hcal_\gamma^{-1}\vec x-1
<
\frac{\Delta_i^2}{4},\quad \Rightarrow\quad |\vec w_i^\conjT\vec x|^2>\frac{3-\Delta_i}{4}.
$$
Since $\Delta_i<1$, the last inequality gives $|\vec w_i^\conjT\vec x|^2>1/2$.
Combining this bound with \eqref{display:nep:product_bound} proves
\eqref{display:nep:wix_estimate}.
\end{proof}

\begin{prop}
\label{prop:nep:left_normalization}
If $\nu_i$ is simple and \eqref{display:nep:assump_gap} holds, the
normalized left eigenvector satisfies
\begin{align}
\label{display:nep:left_norm_bound}
\frac{1}{\normtwo{\vec v_i}^2}
\geq 1-\left|\frac{\impart\eta}{2\Delta_i\re\eta}\right|^2.
\end{align}
For any $\vec x,\vec y\in\{\vec u_i,\vec v_i/\normtwo{\vec v_i}\}$,
the bilinear estimate is
\begin{align}
\label{display:nep:bilinear_residual}
\left|\vec x^\conjT(\Hcal_\gamma-\nu_iI)\vec y\right|
\leq\frac{|\impart\eta|^2q_i}{2(\re\eta)\normtwo{\vec v_i}}
\left(1+\frac{1}{\Delta_i}\right).
\end{align}
\end{prop}
\begin{proof}
We set $p:=\min\{|\vec w_i^\conjT\vec u_i|^2,
|\vec w_i^\conjT\vec v_i|^2/\normtwo{\vec v_i}^2\}$.
By \eqref{display:nep:wix_estimate}, we have $p>1/2$.
In the inner product of the two unit eigenvectors, the term from
$\vec w_i$ has modulus at least $p$. The sum of the other terms has
modulus at most $1-p$ by Cauchy--Schwarz. Since $\vec v_i^\conjT\vec u_i=1$,
the reverse triangle inequality gives
$$
\frac{1}{\normtwo{\vec v_i}}\geq p-(1-p)=2p-1.
$$
For the unit eigenvector attaining $p$, the bounds
\eqref{display:nep:defect_upper} and \eqref{display:nep:defect_lower}
give
$$
4p(1-p)\leq\left|\impart\eta/(2\Delta_i\re\eta)\right|^2.
$$
Thus the normalization bound \eqref{display:nep:left_norm_bound}
follows from
$$
\frac{1}{\normtwo{\vec v_i}^2}\geq (2p-1)^2=1-4p(1-p).
$$

The weighted spectral estimate uses the first inequality in the defect
estimate preceding \eqref{display:nep:defect_upper}. For either unit
eigenvector $\vec x$, this inequality and $|\mu_i|^2\geq(\re\mu_i)^2$ give
\begin{align}
\label{display:nep:defect_q_upper}
(\re\mu_i)\vec x^\conjT\Hcal_\gamma^{-1}\vec x-1
\leq
\frac{|\impart\eta|^2q_i}{(\re\eta)(\re\mu_i)}
\left(1-\frac{(\re\eta)q_i}{\re\mu_i}\right).
\end{align}
By \eqref{display:nep:assump_gap} and \eqref{display:nep:qi_prop},
the ratio $|\impart\eta|/(\Delta_i\re\eta)$ and the factor
$1-(\re\eta)q_i/\re\mu_i$ are at most one. Hence
\eqref{display:nep:defect_q_upper} and \eqref{display:nep:defect_lower},
applied to the unit eigenvector attaining $p$, give
$4p(1-p)\leq(\re\eta)q_i/\re\mu_i$. Since $2p-1\leq p^2$ and
$2p-1\leq 1/\normtwo{\vec v_i}$, we obtain
\begin{align}
\label{display:nep:q_reduction}
1-\frac{(\re\eta)q_i}{\re\mu_i}
\leq (2p-1)^2\leq\frac{p^2}{\normtwo{\vec v_i}}.
\end{align}

\begin{samepage}
For $j\ne i$, the gap condition \eqref{display:nep:assump_gap} gives
$\nu_i+\nu_j\leq|\nu_j-\nu_i|/\Delta_i$. Since
$2\nu_j\leq\nu_i+\nu_j+|\nu_j-\nu_i|$, we have
$$
\left|\nu_j-\nu_i\right|
\leq \frac12\left(1+\frac{1}{\Delta_i}\right)
\frac{(\nu_j-\nu_i)^2}{\nu_j}.
$$
\end{samepage}
The last inequality in the spectral expansion preceding
\eqref{display:nep:defect_lower} bounds the sum with weights
$(\nu_j-\nu_i)^2/\nu_j$. This bound, the weight comparison above, and
$\nu_i|\vec w_i^\conjT\vec x|^2\leq\re\mu_i$ give
\begin{align}
\label{display:nep:weighted_sum_final}
\begin{aligned}
\sum_{j\ne i}|\nu_j-\nu_i|\,|\vec w_j^\conjT\vec x|^2
&\leq\frac12\left(1+\frac{1}{\Delta_i}\right)
\frac{\nu_i}{|\vec w_i^\conjT\vec x|^2}
\left((\re\mu_i)\vec x^\conjT\Hcal_\gamma^{-1}\vec x-1\right)\\
&\leq\frac12\left(1+\frac{1}{\Delta_i}\right)\frac{\re\mu_i}{p^2}
\left((\re\mu_i)\vec x^\conjT\Hcal_\gamma^{-1}\vec x-1\right)\\
&\leq\frac{|\impart\eta|^2q_i}{2(\re\eta)\normtwo{\vec v_i}}
\left(1+\frac{1}{\Delta_i}\right).
\end{aligned}
\end{align}
The last step uses \eqref{display:nep:defect_q_upper} and
\eqref{display:nep:q_reduction}. In the eigenbasis of $\Hcal_\gamma$,
Cauchy--Schwarz with weights $|\nu_j-\nu_i|$ bounds the modulus of the
bilinear form by the geometric mean of the sums in
\eqref{display:nep:weighted_sum_final} for $\vec x$ and $\vec y$.
This gives \eqref{display:nep:bilinear_residual}.
\end{proof}

\begin{prop}
\label{prop:nep:regularity_gradient}
Under the hypotheses of Proposition~\ref{prop:nep:left_normalization},
the material term in \eqref{display:nep:gradient_defect} satisfies
\begin{align}
\label{display:nep:material_defect_bound}
\left|(\vec v_i-\vec u_i)^\conjT\Qcal\vec u_i\right|
\leq\frac{|\impart\eta|q_i}{2\re\eta}
\left(1+\frac{1}{\Delta_i}\right)
\left(1+\frac{1}{\normtwo{\vec v_i}}\right).
\end{align}
\end{prop}
\begin{proof}
If $\impart\eta=0$, simplicity and normalization give $\vec v_i=\vec u_i$.
Otherwise, the right eigenvalue equation and
$(\vec v_i-\vec u_i)^\conjT\vec u_i=0$ give
$(\vec v_i-\vec u_i)^\conjT(\Kcal_\gamma-\nu_iI)\vec u_i=0$.
The splitting $\Kcal_\gamma=\Hcal_\gamma+\imagunit(\impart\eta)\Qcal$
then gives the equality below. The triangle inequality and
\eqref{display:nep:bilinear_residual} with $\vec y=\vec u_i$ and
$\vec x=\vec v_i/\normtwo{\vec v_i}$ or $\vec u_i$ give the bound.
\begin{align*}
|\impart\eta|\,\left|(\vec v_i-\vec u_i)^\conjT\Qcal\vec u_i\right|
&=\left|(\vec v_i-\vec u_i)^\conjT(\Hcal_\gamma-\nu_iI)\vec u_i\right|\\
&\leq(\normtwo{\vec v_i}+1)
\frac{|\impart\eta|^2q_i}{2(\re\eta)\normtwo{\vec v_i}}
\left(1+\frac{1}{\Delta_i}\right).
\end{align*}
The estimate \eqref{display:nep:material_defect_bound} follows upon division
by $|\impart\eta|$.
\end{proof}

\begin{rmk}
\label{rmk:nep:paper3_comparison}
The bound \eqref{display:nep:material_defect_bound} and the identity
\eqref{display:nep:gradient_defect} give
\begin{align}
\label{display:nep:q_grad_bound}
|\nabla q_i|\leq\frac{|\der\eta|q_i}{2\re\eta}
\left(1+\frac1{\Delta_i}\right)
\left(1+\frac1{\normtwo{\vec v_i}}\right).
\end{align}
The strict gap persists locally. If $\eta$ is locally nonconstant,
the open mapping theorem implies that points with $\impart\eta\ne0$
are dense nearby. The local real analyticity of the eigenvectors
therefore extends the bound to $\impart\eta=0$ by continuity.
If $\eta$ is locally constant, then $q_i$ is locally constant.
In the Hermitian setting of \cite{paper3}, $\vec v_i=\vec u_i$ gives
$\der\mu_i=q_i\der\eta$, so \eqref{display:nep:qnewton} is the exact
Newton update and $\rho_i=0$.
Here the left--right derivative defect produces a linear error factor.
For $\re\omega>0$ and $\re\der\eta<0$,
\eqref{display:nep:rho_grad_abs}, \eqref{display:nep:denominator_bound}, and
\eqref{display:nep:q_grad_bound} give
\begin{align}
\label{display:nep:rho_relative_bound}
|\rho_i|\leq
\frac{|\impart\eta|\,|\der\eta|q_i}
{2\re\eta\left(2\re\omega-(\re\der\eta)q_i\right)}
\left(1+\frac1{\Delta_i}\right)
\left(1+\frac1{\normtwo{\vec v_i}}\right).
\end{align}
For exact inner eigenpairs, a value of the right-hand side below one
at a simple root guarantees local contraction of the update
\eqref{display:nep:qnewton} by \eqref{display:nep:err_exp_complete}.
\end{rmk}

\section{Numerical Experiments}
\label{se:num}

\subsection{Parameters and Test Examples}
\label{sse:num:parameters}

We consider the simple cubic (SC), face centered cubic (FCC),
and body centered cubic (BCC) lattices as in \cite{paper2}.
We use dimensionless coordinates and frequencies, with the conventional
cubic lattice constant and the background wave speed normalized to one.
The primitive
lattice vectors are defined by
\begin{align}
\label{display:num:lattice}
\begin{aligned}
\text{SC lattice:}&\quad\vec a_1=(1,0,0)^\T,\quad
 \vec a_2=(0,1,0)^\T,\quad
 \vec a_3=(0,0,1)^\T,\\
\text{FCC lattice:}&\quad\vec a_1=\left(0,\tfrac12,\tfrac12\right)^\T,\quad
 \vec a_2=\left(\tfrac12,0,\tfrac12\right)^\T,\quad
 \vec a_3=\left(\tfrac12,\tfrac12,0\right)^\T,\\
\text{BCC lattice:}&\quad\vec a_1=\left(-\tfrac12,\tfrac12,\tfrac12\right)^\T,\quad
 \vec a_2=\left(\tfrac12,-\tfrac12,\tfrac12\right)^\T,\quad
 \vec a_3=\left(\tfrac12,\tfrac12,-\tfrac12\right)^\T.
\end{aligned}
\end{align}
Their high-symmetry points are
\begin{align}\label{display:num:highsym}
\begin{aligned}
\text{SC lattice:}&\quad
 K(0,0,0),\quad L(\pi,0,0),\quad M(\pi,\pi,0),\quad N(\pi,\pi,\pi),\\
\text{FCC lattice:}&\quad
 X(0,2\pi,0),\quad U\left(\tfrac{\pi}{2},2\pi,\tfrac{\pi}{2}\right),\quad
 L(\pi,\pi,\pi),\\
&\quad\Gamma(0,0,0),\quad W(\pi,2\pi,0),\quad
 K\left(\tfrac{3\pi}{2},\tfrac{3\pi}{2},0\right),\\
\text{BCC lattice:}&\quad H(0,0,2\pi),\quad\Gamma(0,0,0),\quad
 P(\pi,\pi,\pi),\quad N(\pi,0,\pi).
\end{aligned}
\end{align}

The wave number vector follows
$K\longrightarrow L\longrightarrow M\longrightarrow N\longrightarrow K$
for the SC lattice and
$X\longrightarrow U\longrightarrow L\longrightarrow\Gamma
\longrightarrow X\longrightarrow W\longrightarrow K$
for the FCC lattice. The BCC path is
$H\longrightarrow\Gamma\longrightarrow P\longrightarrow H
\longrightarrow N\longrightarrow\Gamma\longrightarrow H
\longrightarrow P\longrightarrow N$.
The first two occurrences of $H$ denote $(0,0,2\pi)$, and the last
denotes the symmetry-equivalent point $(0,2\pi,0)$.
Each path segment contains 20 equally spaced samples, including
its first endpoint and excluding its last endpoint. Therefore,
one complete band calculation contains 80 wave number vectors for the
SC lattice, 120 for the FCC lattice, and 160 for the BCC lattice.

We consider five material regions. The region
\texttt{SC\_CURV} consists of a central sphere of radius $0.345$ and
three cylinders of radius $0.11$ \cite{paper2}:
\begin{align}
\label{display:num:sc_curv_region}
\begin{aligned}
\Omega_{\texttt{SC\_CURV}}
&=\left\{\left(x-\tfrac12\right)^2+\left(y-\tfrac12\right)^2
+\left(z-\tfrac12\right)^2\leq0.345^2\right\}\\
&\quad\cup\left\{\left(x-\tfrac12\right)^2+\left(y-\tfrac12\right)^2\leq0.11^2\right\}
\cup\left\{\left(y-\tfrac12\right)^2+\left(z-\tfrac12\right)^2\leq0.11^2\right\}\\
&\quad\cup\left\{\left(z-\tfrac12\right)^2+\left(x-\tfrac12\right)^2\leq0.11^2\right\}.
\end{aligned}
\end{align}
 
The region
\texttt{SC\_FLAT} is the following union of four cuboids:
\begin{align}
\label{display:num:sc_flat_region}
\begin{aligned}
\Omega_{\texttt{SC\_FLAT}}
&=
[0,\tfrac14]\times[0,\tfrac14]\times[0,1]
\ \cup\
[0,\tfrac14]\times[0,1]\times[\tfrac14,\tfrac12]\\
&\quad\cup\
[0,1]\times[\tfrac12,\tfrac34]\times[\tfrac12,\tfrac34]
\ \cup\
[\tfrac12,\tfrac34]\times[0,1]\times[\tfrac34,1].
\end{aligned}
\end{align}
The region \texttt{FCC} adopts the diamond-type configuration used in
\cite{huang2013eigendecomppc,paper2}. It consists of spheres of radius
$0.12$ and connecting spheroids with semi-minor axis $0.11$. The foci of each spheroid are located
at the centers of two neighboring spheres.
The sphere centers are the FCC lattice sites and their translates
by $\left(\tfrac14,\tfrac14,\tfrac14\right)^\T$, with periodic copies included.

The \texttt{BCC-SG} and \texttt{BCC-DG} regions use the single and
double gyroid geometries in \cite{lu2013weyl,paper2}.
The double gyroid region is defined by
\begin{align}
\label{display:num:bcc_dg_region}
\begin{aligned}
\Omega_{\texttt{BCC-DG}}
&=\left\{(x,y,z):
\left|\sin(2\pi x)\cos(2\pi y)+\sin(2\pi y)\cos(2\pi z)
\right.\right.\\
&\hspace{7.5em}\left.\left.+\sin(2\pi z)\cos(2\pi x)\right|>1.1\right\}.
\end{aligned}
\end{align}
The single gyroid \texttt{BCC-SG} retains the positive network,
where the trigonometric sum in \eqref{display:num:bcc_dg_region}
is greater than $1.1$.
The \texttt{FCC}, \texttt{BCC-SG}, and \texttt{BCC-DG} examples use
\texttt{POL}.

The implementation applies the inverse dielectric coefficient
$\eta(\omega)=\varepsilon^{-1}(\omega)$.
The \texttt{SCAL} and \texttt{SiC} models use the scalar Lorentz
law \cite{effenberger2012absorbing} and the ionic Lorentz law for
silicon carbide \cite{devarapu2012sic}, respectively.
The \texttt{POL} and \texttt{MIX} models are artificial rational test
functions. The \texttt{POL} form is motivated by the lossless GaAs
phonon-polariton law \cite{sigalas1994dispersive}; its complex
quadratic coefficient is artificial.
Table~\ref{tab:num:dielectric_models} lists the functions and parameters.
The material parameters $\gamma_{\mathrm d}$ and $\Gamma$ in this table
are independent of the kernel-compensation coefficient $\gamma$.

\begin{table}[pos={!htbp}]
\centering
\caption{The dielectric functions with normalized parameters used in our numerical examples.}
\label{tab:num:dielectric_models}
\normalfont\fontsize{11}{13.6}\selectfont
\renewcommand{\small}{\fontsize{11}{13.6}\selectfont}

\renewcommand{\arraystretch}{1.3}
\setlength{\extrarowheight}{2pt}

\small
\setlength{\tabcolsep}{4pt}
\begin{tabular*}{\linewidth}{
  @{\extracolsep{\fill}}c
     >{\raggedright\arraybackslash}p{0.42\textwidth}
     >{\raggedright\arraybackslash}p{0.36\textwidth}@{}}
\toprule
Model
&
Function $\varepsilon(\omega)$
&
Parameters
\\
\midrule

\texttt{SCAL}
&
$\displaystyle
\varepsilon(\omega)
=\varepsilon_\infty+
\frac{S}
{\omega_0^2-\omega^2+\imagunit\gamma_{\mathrm d}\omega}
$
&
$(\varepsilon_\infty,\omega_0,S,\gamma_{\mathrm d})=(6,5.60,12,0.06)$
\\

\texttt{SiC}
&
$\displaystyle
\varepsilon(\omega)
=\varepsilon_\infty\left(
1+
\frac{\omega_p^2}
{\omega_0^2-\omega^2+\imagunit\Gamma\omega}
\right)$
& $(\varepsilon_\infty,\omega_0,\omega_p,\Gamma)=(6.70, 5.13, 3.60, 0.03)$
\\

\texttt{POL}
&
$\displaystyle
\varepsilon(\omega)
=\varepsilon_\infty
\frac{\omega_{\mathrm L}^2-\omega^2}
{\omega_{\mathrm T}^2-(1+\imagunit\gamma_{\mathrm d})\omega^2}$
&
$(\varepsilon_\infty,\omega_{\mathrm T},\omega_{\mathrm L},\gamma_{\mathrm d})=(20,8.12,8.75,0.02)$
\\

\texttt{MIX}
&
$\displaystyle
\varepsilon(\omega)
=11.8+
\frac{22.8-1.2\omega+2\omega^2+0.11\omega^3}
{6.7-0.3\omega-0.025\omega^2+\imagunit\gamma_{\mathrm d}\omega}$
&
$\gamma_{\mathrm d}=0.01$
\\

\bottomrule
\end{tabular*}
\end{table}

The \texttt{SiC} parameters in Table~\ref{tab:num:dielectric_models}
are rounded. The computations use $\omega_0=2\pi(0.817)$,
$\omega_p=2\pi\sqrt{0.998^2-0.817^2}$, and
$\Gamma=2\pi(0.004904)$.
The Lorentz damping terms follow the
$\exp(+\imagunit\omega t)$ convention.
Since $\left|\exp(+\imagunit\omega t)\right|=\exp(-\impart\omega\,t)$,
frequencies with positive imaginary part correspond to temporal decay.
The negative imaginary parts in the \texttt{POL} results reflect its
artificial material law.

The band and tolerance experiments use $N=120$ and compute the first
$\nev=10$ eigenpairs at each wave number vector. The Generalized Davidson
method limits the search dimension to four times the retained block size
and allows at most 2000 inner iterations,
except for the \texttt{SC\_CURV}/\texttt{SCAL} band run and the tolerance
experiments, which allow 500. The quasi-Newton iteration allows at most
1000 outer iterations. The complete bands use outer and inner residual
tolerances of $10^{-5}$ and $9\times10^{-6}$, respectively. In the tolerance
experiments, the inner tolerance is $0.9$ times the outer tolerance.

Each run starts with a random block. At each subsequent Bloch wave vector,
the solver reuses the previously computed vectors and fills the remaining
columns with random vectors.
The \texttt{SiC} and \texttt{SC\_CURV}/\texttt{MIX} band calculations start
with 13 vectors. The other band calculations start with four vectors and
enlarge the block to 13 as additional frequencies are computed.
The first frequency is initialized at $\omega_1^{(0)}=0$, so its first
update is a fixed-point step. The implementations used for the convergence
histories and the $\rho_i^*$ diagnostics additionally halve an inadmissible
fallback step toward the preceding frequency, up to 32 times, and retain
the preceding frequency if it is admissible and these trials fail.

The band and tolerance runtimes are obtained with the C++ CUDA implementation
in complex double precision, using CUDA 12.4, cuBLAS, and cuFFT on one
NVIDIA GeForce RTX 4090D GPU with 24 GiB of device memory. The timing
experiments do not perform the optional adjoint-eigenvector postprocessing.

\subsection{Figures of Complete Band}
\label{sse:num:complete_bands}

The upper and lower panels show $\frac{|\re\omega|}{2\pi}$ and
$\frac{\impart\omega}{2\pi}$, respectively. Unconnected markers are used
because no mode matching is applied between adjacent wave number vectors.

A sampled bandgap is determined from the absolute real parts of the
computed frequencies along each path. Its ratio is defined by
\begin{align}
\label{display:num:bandgap_ratio}
\text{bandgap ratio}
=2\frac{\omega_{\mathrm{up}}-\omega_{\mathrm{low}}}
{\omega_{\mathrm{up}}+\omega_{\mathrm{low}}}.
\end{align}
Here, $\omega_{\mathrm{low}}$ is the largest frequency below the gap,
and $\omega_{\mathrm{up}}$ is the smallest frequency above the gap.
Shading marks gaps along the sampled paths; it does not establish
spectral exclusion throughout the Brillouin zone.
The detection criterion uses the 95th percentile of changes between
adjacent wave number samples, after sorting the absolute real parts
at each sample, as an estimate of the sampling resolution.
We retain intervals whose lower edge is at least one quarter of the
largest displayed frequency, whose width exceeds $1.5$ times this
resolution, and whose ratio is at least $0.05$.
The widest retained interval is reported.

\begin{figure}[pos={!htbp}]
\centering
\begin{subfigure}[t]{0.48\textwidth}
\centering
\includegraphics[width=\linewidth,height=0.22\textheight,keepaspectratio]
{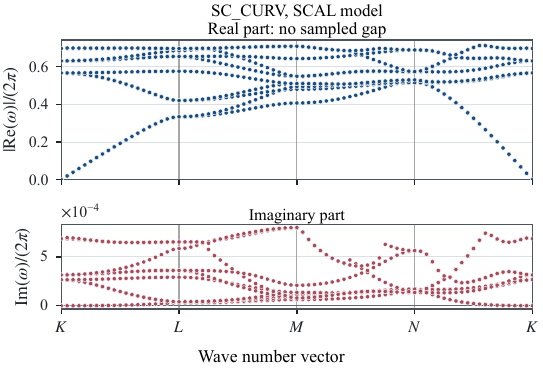}
\caption{The \texttt{SC\_CURV}/\texttt{SCAL} example has no detected sampled bandgap.}
\label{fig:num:scal}
\end{subfigure}
\hfill
\begin{subfigure}[t]{0.48\textwidth}
\centering
\includegraphics[width=\linewidth,height=0.22\textheight,keepaspectratio]
{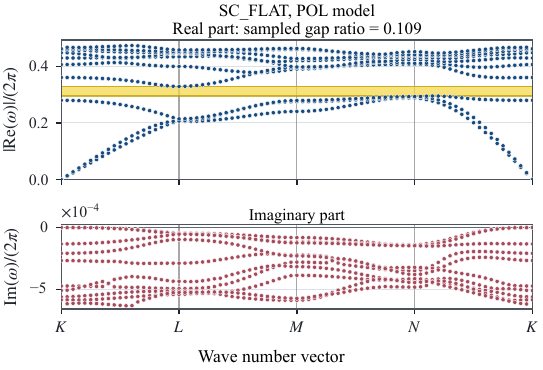}
\caption{The \texttt{SC\_FLAT}/\texttt{POL} example has bandgap ratio $0.109$.}
\label{fig:num:polar}
\end{subfigure}
\caption{Complete complex bands for the two SC examples.}
\label{fig:num:sc_examples}
\end{figure}

\begin{figure}[pos={!htbp}]
\centering
\begin{subfigure}[t]{0.48\textwidth}
\centering
\includegraphics[width=\linewidth,height=0.22\textheight,keepaspectratio]
{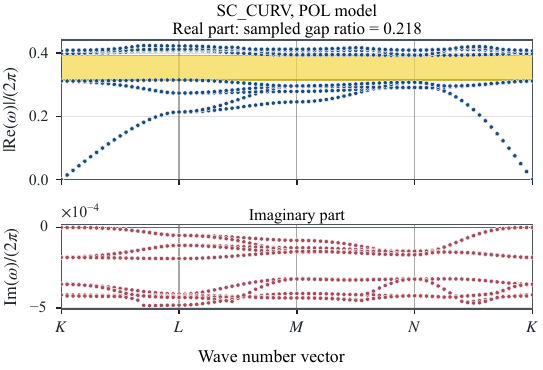}
\caption{The \texttt{POL} example has sampled bandgap ratio $0.218$.}
\label{fig:num:sc_curv_pol}
\end{subfigure}
\hfill
\begin{subfigure}[t]{0.48\textwidth}
\centering
\includegraphics[width=\linewidth,height=0.22\textheight,keepaspectratio]
{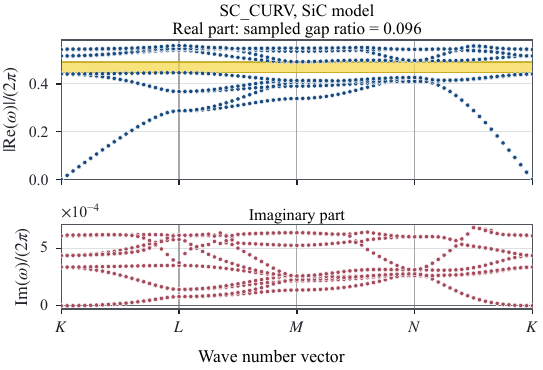}
\caption{The \texttt{SiC} example has sampled bandgap ratio $0.096$.}
\label{fig:num:sic}
\end{subfigure}
\caption{Complete complex bands for the \texttt{SC\_CURV}/\texttt{POL} and \texttt{SC\_CURV}/\texttt{SiC} examples.}
\label{fig:num:sc_curv_pol_sic}
\end{figure}

Figures~\ref{fig:num:sc_examples} and~\ref{fig:num:sc_curv_pol_sic}
compare the SC examples. For the fixed \texttt{SC\_CURV} geometry,
\texttt{SCAL} has no detected sampled gap, while \texttt{POL} and
\texttt{SiC} give ratios of $0.218$ and $0.096$, respectively.
For \texttt{POL}, changing the material region from
\texttt{SC\_CURV} to \texttt{SC\_FLAT} reduces the ratio to $0.109$.
These comparisons show the dependence of the sampled gaps on both
the material law and the material geometry.

\begin{figure}[pos={!htbp}]
\centering
\captionsetup{justification=centering}
\begin{subfigure}[t]{0.48\textwidth}
\centering
\includegraphics[width=\linewidth,height=0.22\textheight,keepaspectratio]
{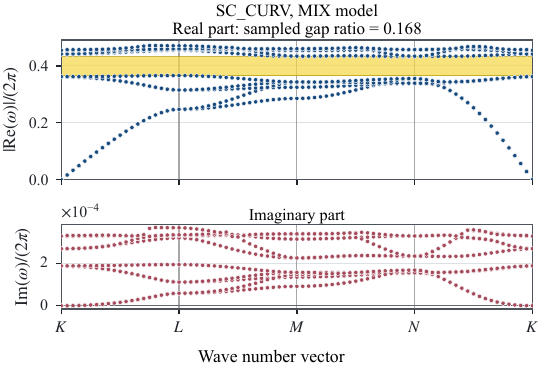}
\caption{The \texttt{SC\_CURV}/\texttt{MIX} bands have a sampled
bandgap ratio of $0.168$.}
\label{fig:num:sc_curv_mix}
\end{subfigure}
\hfill
\begin{subfigure}[t]{0.48\textwidth}
\centering
\includegraphics[width=\linewidth,height=0.22\textheight,keepaspectratio]
{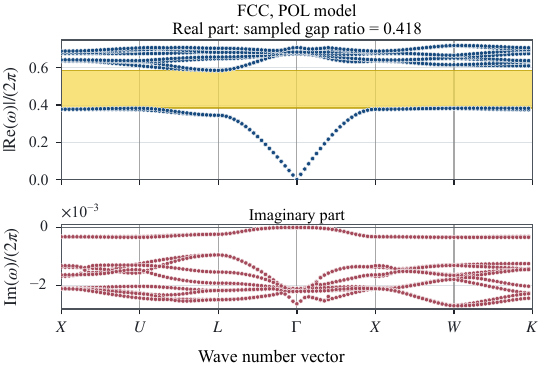}
\caption{The \texttt{FCC}/\texttt{POL} bands have a sampled bandgap
ratio of $0.418$.}
\label{fig:num:fcc_pol}
\end{subfigure}
\caption{The bands for the \texttt{SC\_CURV}/\texttt{MIX} and
\texttt{FCC}/\texttt{POL} examples are shown.}
\label{fig:num:sc_curv_mix_fcc_pol}
\end{figure}

Figure~\ref{fig:num:sc_curv_mix_fcc_pol} shows a sampled gap ratio of
$0.168$ for \texttt{SC\_CURV}/\texttt{MIX} and $0.418$ for
\texttt{FCC}/\texttt{POL}. The latter is the largest ratio among the
ten examples.

\begin{figure}[pos={!htbp}]
\centering
\begin{subfigure}[t]{0.48\textwidth}
\centering
\includegraphics[width=\linewidth,height=0.22\textheight,keepaspectratio]
{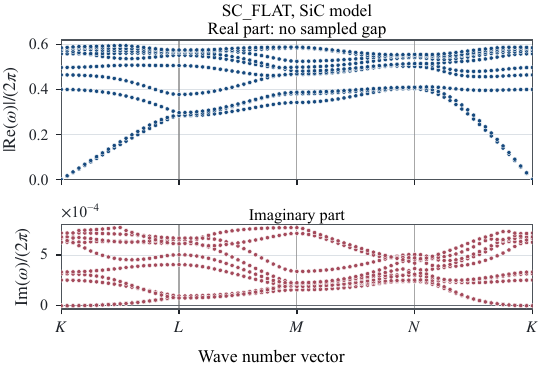}
\caption{The \texttt{SiC} example has no detected sampled bandgap.}
\label{fig:num:sc_flat_sic}
\end{subfigure}
\hfill
\begin{subfigure}[t]{0.48\textwidth}
\centering
\includegraphics[width=\linewidth,height=0.22\textheight,keepaspectratio]
{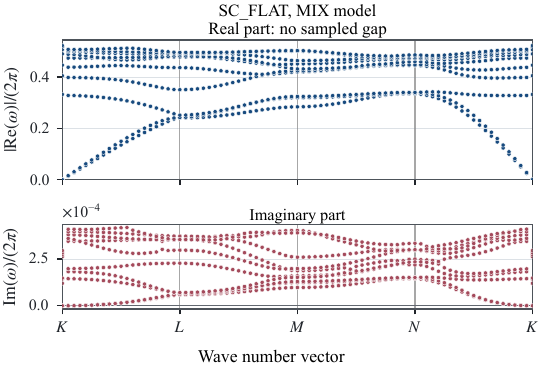}
\caption{The \texttt{MIX} example has no detected sampled bandgap.}
\label{fig:num:sc_flat_mix}
\end{subfigure}
\caption{Complete complex bands for the \texttt{SC\_FLAT}/\texttt{SiC} and \texttt{SC\_FLAT}/\texttt{MIX} examples.}
\label{fig:num:sc_flat_sic_mix}
\end{figure}

Figure~\ref{fig:num:sc_flat_sic_mix} shows no detected sampled gap for
either \texttt{SiC} or \texttt{MIX} on \texttt{SC\_FLAT}.
Both laws give gaps on \texttt{SC\_CURV}, so a change of geometry alone
changes whether a gap is detected along the corresponding sampled path.

\begin{figure}[pos={!htbp}]
\centering
\begin{subfigure}[t]{0.48\textwidth}
\centering
\includegraphics[width=\linewidth,height=0.22\textheight,keepaspectratio]
{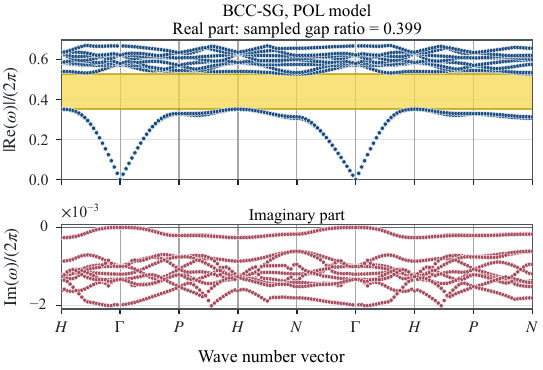}
\caption{The sampled bandgap ratio is $0.399$.}
\label{fig:num:bcc_sg_pol}
\end{subfigure}
\hfill
\begin{subfigure}[t]{0.48\textwidth}
\centering
\includegraphics[width=\linewidth,height=0.22\textheight,keepaspectratio]
{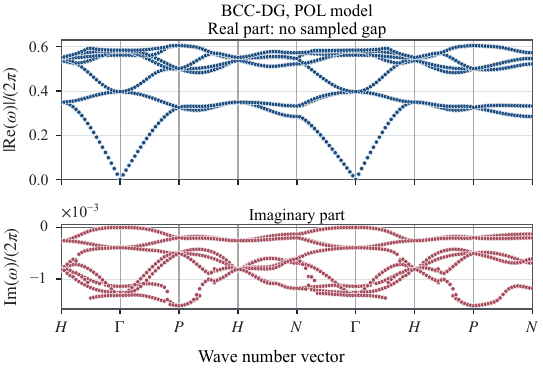}
\caption{No sampled bandgap is detected.}
\label{fig:num:bcc_dg_pol}
\end{subfigure}
\caption{The \texttt{BCC-SG}/\texttt{POL} and
\texttt{BCC-DG}/\texttt{POL} bands are compared.}
\label{fig:num:bcc_examples}
\end{figure}

Figure~\ref{fig:num:bcc_examples} compares two BCC material regions
with the same \texttt{POL} law. The single gyroid has a sampled gap
ratio of $0.399$, whereas no sampled gap is detected for the double
gyroid. The different outcomes on the same lattice further illustrate
the effect of the material geometry.

Figures~\ref{fig:num:sc_examples}--\ref{fig:num:bcc_examples}
show the ten complex band diagrams.
Table~\ref{tab:num:runtime} summarizes the runtimes of the ten
complete calculations.
\FloatBarrier

\subsection{Runtime and Performance}
\label{sse:num:runtime}

Runtime includes initialization and the quasi-Newton solve for the
first ten frequencies, with GPU synchronization; file output is excluded. Common initialization times are summed
and divided equally among all sampled wave number vectors.
Table~\ref{tab:num:runtime} reports the mean and population standard
deviation along each path.

\begin{table}[pos={!htbp}]
\centering
\caption{Mean and population standard deviation of runtime (s) per
wave number vector at tolerance $10^{-5}$, including initialization.}
\label{tab:num:runtime}
\normalfont\fontsize{11}{13.6}\selectfont
\renewcommand{\small}{\fontsize{11}{13.6}\selectfont}

\renewcommand{\arraystretch}{1.3}
\setlength{\extrarowheight}{2pt}

\small
\setlength{\tabcolsep}{2pt}
\begin{tabular*}{\linewidth}{@{\extracolsep{\fill}}l*{10}{c}@{}}
\toprule
Example & \multicolumn{4}{c}{\texttt{SC\_CURV}} & \multicolumn{3}{c}{\texttt{SC\_FLAT}} & \texttt{FCC} & \texttt{BCC-SG} & \texttt{BCC-DG} \\
\cmidrule(lr){2-5}\cmidrule(lr){6-8}
Model & \texttt{SCAL} & \texttt{SiC} & \texttt{POL} & \texttt{MIX} & \texttt{SiC} & \texttt{POL} & \texttt{MIX} & \texttt{POL} & \texttt{POL} & \texttt{POL} \\
\midrule
Mean (s) & $109.700$ & $215.489$ & $150.076$ & $241.571$ & $283.963$ & $190.815$ & $315.947$ & $351.551$ & $349.752$ & $249.936$ \\
Std (s) & $17.807$ & $47.884$ & $25.279$ & $30.709$ & $20.689$ & $11.854$ & $23.736$ & $69.539$ & $37.627$ & $49.044$ \\
\bottomrule
\end{tabular*}
\end{table}

We use the \texttt{SC\_CURV}/\texttt{SCAL} example at
$\wnk=\left(\frac{\pi}{2},\pi,\frac{\pi}{5}\right)$ to compare the outer tolerances
$10^{-4}$, $10^{-5}$, $10^{-6}$, and $10^{-7}$.
Table~\ref{tab:num:tolerance_runtime} reports the runtimes and the
maximum absolute differences in the first ten normalized complex
frequencies $\frac{\omega}{2\pi}$ from the matched $10^{-7}$ results.
The runtimes include initialization and are taken from the first run.

\begin{table}[pos={!htbp}]
\centering
\caption{Runtime and maximum normalized frequency difference
from the result at tolerance \texttt{1e-7}.}
\label{tab:num:tolerance_runtime}
\normalfont\fontsize{11}{13.6}\selectfont
\renewcommand{\small}{\fontsize{11}{13.6}\selectfont}

\renewcommand{\arraystretch}{1.3}
\setlength{\extrarowheight}{2pt}

\begin{tabular*}{\linewidth}{@{\extracolsep{\fill}}lcccc@{}}
\toprule
Tolerance & \texttt{1e-4} & \texttt{1e-5} & \texttt{1e-6} & \texttt{1e-7} \\
\midrule
Runtime (s) & $104.465$ & $143.343$ & $196.719$ & $252.089$ \\
Max diff from tol \texttt{1e-7} & \texttt{7.778e-11} & \texttt{1.771e-11}
& \texttt{2.192e-12} & $0$ \\
\bottomrule
\end{tabular*}
\end{table}

The maximum difference over all tolerances is $7.778\times10^{-11}$.
Thus, $10^{-4}$ is sufficient for the frequency accuracy observed
in this fixed-grid example.

We compare the CPU and GPU implementations for the
\texttt{SC\_CURV}/\texttt{SiC} example at $\wnk=(\pi,0,0)$,
with ten frequencies and outer tolerance $10^{-5}$.
The CPU implementation uses one Intel Xeon Gold 6530 core and one
thread with serial FFTW 3.3.10 and Eigen 3.4. The GPU implementation
uses one RTX 4090D with cuFFT and cuBLAS.
Both implementations use complex double precision and the same
random initial block. The Generalized Davidson method limits the search
dimension to four times the retained block size for $N=64,80,96,112$
and three times the retained block size for $N=128$.
Each configuration is timed once. Table~\ref{tab:num:cpu_gpu}
reports the initialization and quasi-Newton solve times.

\begin{table}[pos={!htbp}]
\centering
\caption{CPU and GPU runtimes including initialization for the
\texttt{SC\_CURV}/\texttt{SiC} example.}
\label{tab:num:cpu_gpu}
\normalfont\fontsize{11}{13.6}\selectfont
\renewcommand{\small}{\fontsize{11}{13.6}\selectfont}

\renewcommand{\arraystretch}{1.3}
\setlength{\extrarowheight}{2pt}

\begin{tabular*}{\linewidth}{@{\extracolsep{\fill}}lccccc@{}}
\toprule
$N$ & 64 & 80 & 96 & 112 & 128 \\
\midrule
CPU (s) & $3465.054$ & $6368.870$ & $10950.883$ & $17081.481$ & $28450.894$ \\
GPU (s) & $47.461$ & $91.507$ & $157.164$ & $242.990$ & $405.929$ \\
Speedup & $73.009$ & $69.599$ & $69.678$ & $70.297$ & $70.088$ \\
\bottomrule
\end{tabular*}
\end{table}

We fit the runtime in seconds against the number of degrees of freedom, $3N^3$, by linear least squares.
The dimensionless coefficients of determination are $R^2=0.9934$ for the CPU and $R^2=0.9935$ for the GPU.
These values are close to one and indicate approximately linear growth
over the five tested grid sizes, $64\leq N\leq128$.
The fitted slopes are $4.4894\times10^{-3}$ and
$6.4175\times10^{-5}$ seconds per degree of freedom for the CPU and
GPU, respectively. These fits describe the measured total runtimes over
this range; a single matrix-free block application has the
$O(mN^3\log N)$ cost given in Section~\ref{se:disc}.

\FloatBarrier

\subsection{Convergence Histories}
\label{sse:num:histories}

\begin{figure}[pos={!htbp}]
\centering
\includegraphics[width=0.47\textwidth]{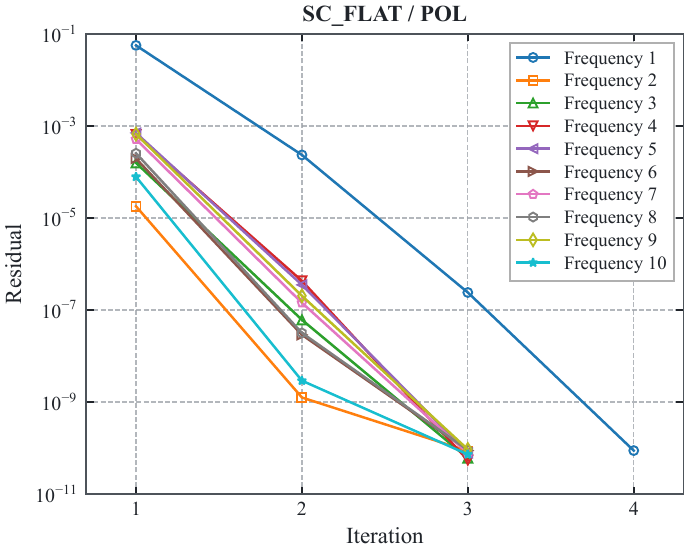}
\hfill
\includegraphics[width=0.47\textwidth]{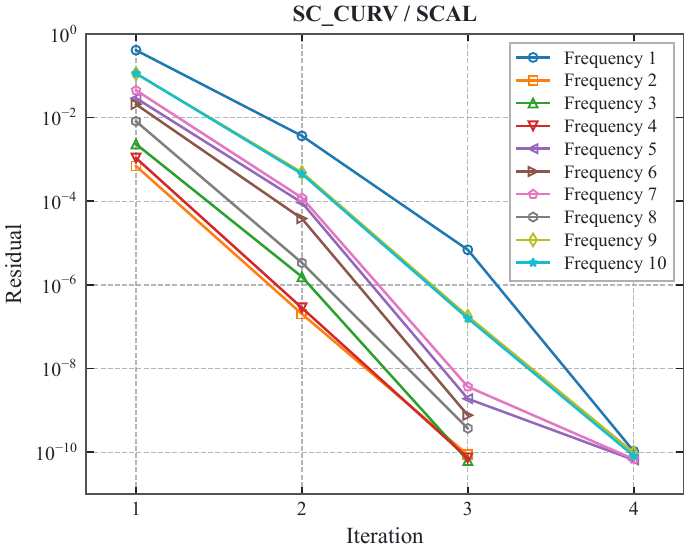}\\
\includegraphics[width=0.47\textwidth]{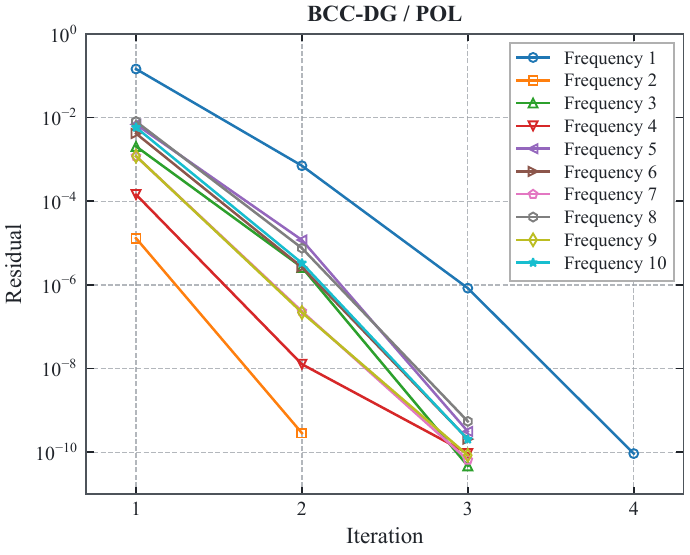}
\hfill
\includegraphics[width=0.47\textwidth]{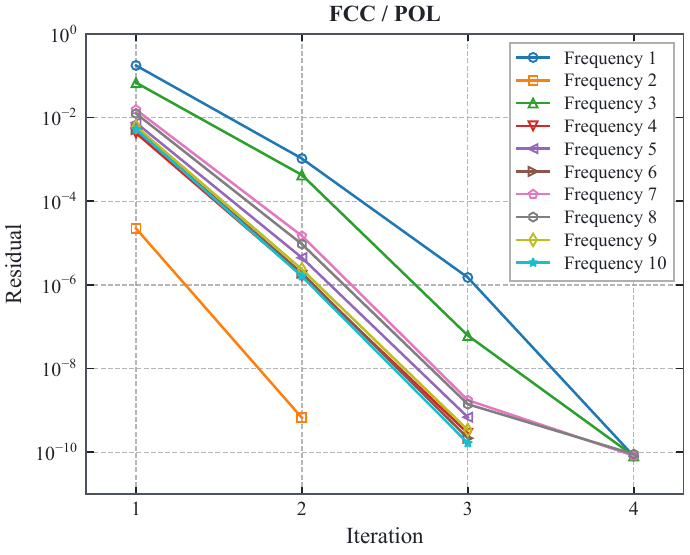}
\caption{Quasi-Newton residual histories for the first ten frequencies.
The examples are \texttt{SC\_FLAT}/\texttt{POL},
\texttt{SC\_CURV}/\texttt{SCAL}, \texttt{BCC-DG}/\texttt{POL},
and \texttt{FCC}/\texttt{POL}, in row order.}
\label{fig:num:residual_histories}
\end{figure}

We compute the first ten frequencies with $N=120$ at
$\wnk=\left(\frac{\pi}{2},\pi,\frac{\pi}{5}\right)$, using outer and inner tolerances
$10^{-9}$ and $10^{-10}$, respectively.
Figure~\ref{fig:num:residual_histories} shows the quasi-Newton
residual histories for four examples.

All forty frequencies converge within two to four iterations,
with final residuals below $7.64\times10^{-10}$.

We evaluate $\rho_i^*$ from \eqref{display:nep:rho} for the
\texttt{SC\_CURV}/\texttt{SCAL} example with $N=120$ and
$\wnk=\left(\frac{\pi}{2},\pi,\frac{\pi}{5}\right)$.
The outer and inner tolerances are $10^{-9}$ and $10^{-10}$, respectively.
The left eigenvectors are computed at the converged frequencies, with
$\normtwo{\vec u_i}=1$ and $\vec v_i^\conjT\vec u_i=1$.
Table~\ref{tab:num:rho} reports the moduli $|\rho_i^*|$ to three significant digits.

\begin{table}[pos={!htbp}]
\centering
\caption{The moduli $|\rho_i^*|$ evaluated from
\eqref{display:nep:rho} for the \texttt{SC\_CURV}/\texttt{SCAL}
example with $N=120$ and $\wnk=\left(\frac{\pi}{2},\pi,\frac{\pi}{5}\right)$.}
\label{tab:num:rho}
\normalfont\fontsize{11}{13.6}\selectfont
\renewcommand{\small}{\fontsize{11}{13.6}\selectfont}

\renewcommand{\arraystretch}{1.3}
\setlength{\extrarowheight}{2pt}

\setlength{\tabcolsep}{2pt}
\begin{tabular*}{\linewidth}{@{\extracolsep{\fill}}c*{10}{c}@{}}
\toprule
$i$ & 1 & 2 & 3 & 4 & 5 & 6 & 7 & 8 & 9 & 10 \\
\midrule
$|\rho_i^*|$ & \text{2.48e-6} & \text{3.73e-6} & \text{1.05e-5}
& \text{1.39e-5} & \text{1.78e-5} & \text{1.95e-5} & \text{2.95e-5}
& \text{1.08e-4} & \text{3.64e-4} & \text{3.50e-4} \\
\bottomrule
\end{tabular*}
\end{table}
\FloatBarrier

\section{Conclusion}
\label{se:con}

We have extended the kernel compensation framework to nonlinear,
non-Hermitian Maxwell eigenproblems with an isotropic dispersive material
law. The fixed discrete curl, divergence, and harmonic spaces give an
exact unitary block diagonalization. Under the stated penalty condition,
this decomposition proves that compensation preserves the eigenpairs of
the original constrained nonlinear problem. The complex-coefficient
preconditioning analysis gives a normal two-sided form, a line-segment
spectral inclusion, and a condition-number bound independent of the grid
size. These results justify the use of the inherited mimetic and DFT
operators beyond the Hermitian setting.

The right Generalized Davidson solver supplies the eigenpair used in each
frequency update. The derivative approximation uses the normalized right
eigenvector and avoids an additional adjoint eigensolve. With exact inner eigenpairs, the iteration is locally linearly
convergent when its error factor has modulus between zero and one. It is at least quadratic when
the derivative defect vanishes at the root. We express this factor through the gradient
of the material Rayleigh quotient and bound it using the relative gap of
the Hermitian part.

The experiments demonstrate the complete matrix-free method for ten
complex band calculations across five material geometries and three
lattice types. In the four convergence tests, all forty frequencies
converge within two to four outer iterations. The measured linear factors
in the reported test are small, consistent with rapid convergence at the
computed accuracy. For the specified \texttt{SC\_CURV}/\texttt{SiC}
comparison, the GPU implementation achieves speedups of $69.6$--$73.0$
over the single-core, single-thread CPU implementation for $64\leq N\leq128$.
These results support the practical use of kernel compensation and the
right-Rayleigh quasi-Newton update for large three-dimensional dispersive
band calculations. An extension of the local analysis to multiple
eigenvalues requires an invariant-subspace formulation.

\section*{Declaration of competing interest}
\label{sec:competing_interest}
The authors declare that they have no known competing financial interests
or personal relationships that could have influenced the work reported
in this paper.

\section*{Data availability}
\label{sec:data_availability}
The data and code supporting this study are available from the corresponding
author upon reasonable request.

\section*{Declaration of generative AI and AI-assisted technologies in the manuscript preparation process}
\label{sec:ai_declaration}
During the preparation of this manuscript, the authors used OpenAI Codex
to assist with manuscript revision, mathematical exposition, notation
checks, and LaTeX preparation. The authors direct and review the revisions
and take full responsibility for the content of the manuscript.

\bibliographystyle{abbrv}
\bibliography{ref}
\end{document}

%% file: preamble.tex
\usepackage{indentfirst}
\usepackage{amsmath}
\usepackage{amssymb}
\usepackage{amsthm}
\usepackage{mathtools}

\numberwithin{equation}{section}

\usepackage{mathrsfs}

\input{ursfs.fd}
\DeclareFontShape{U}{rsfs}{m}{n}{%
  <-6.5>rsfs5 <6.5-8>rsfs7 <8->rsfs10
}{}

\allowdisplaybreaks

\usepackage{graphicx}
\usepackage{placeins}
\graphicspath{{./}}

\usepackage{caption}
\usepackage{subcaption}

\usepackage{array}
\usepackage{multirow}
\usepackage{booktabs}
\usepackage{colortbl}

\usepackage{tikz}
\usepackage{tikz-cd}
\usepackage{tikz-3dplot}

\usepackage[all,pdf]{xy}

\usepackage{pgffor}

\usepackage{algorithm}
\usepackage{algpseudocode}

\usepackage{listings}

\usepackage{setspace}
\usepackage{multicol}
\usepackage{comment}
\usepackage{lipsum}

\usepackage{xcolor}
\usepackage{xurl}

\usepackage{chngcntr}
\usepackage[numbers,sort&compress]{natbib}

\counterwithout{figure}{section}
\counterwithout{algorithm}{section}

\counterwithout{table}{section}

\newtheorem{thm}{Theorem}[section]
\newtheorem{lem}[thm]{Lemma}
\newtheorem{prop}[thm]{Proposition}

\newtheorem{rmk}{Remark}[section]

\foreach \x in {A,...,Z}{
    \expandafter\xdef\csname \x cal\endcsname{
        {\noexpand\mathcal{\x}}
    }
}

\foreach \x in {A,...,Z}{
    \expandafter\xdef\csname \x bf\endcsname{
        {\noexpand\mathbf{\x}}
    }
}

\foreach \x in {A,...,Z}{
    \expandafter\xdef\csname \x bb\endcsname{
        {\noexpand\mathbb{\x}}
    }
}

\renewcommand{\vec}[1]{\boldsymbol{#1}}

\newcommand{\p}{\partial}

\DeclareMathOperator{\diag}{diag}
\DeclareMathOperator{\rank}{rank}

\DeclareMathOperator{\image}{im}

\DeclareMathOperator{\impart}{Im}

\usepackage[normalem]{ulem}

\DeclarePairedDelimiter{\norm}{\lVert}{\rVert}

\newcommand{\normtwo}[1]{
    \norm*{#1}_2
}

\newcommand{\wnk}{\vec{\mathrm{k}}}
\newcommand{\wnkk}{\mathrm{k}}

\newcommand{\conjT}{\dagger}

\newcommand{\T}{\top}

\newcommand{\imagunit}{\mathrm{i}}

\DeclareMathOperator{\re}{Re}

\newcommand{\nev}{\mathrm{nev}}

\usepackage{ifthen}

\newcommand{\useprime}{1}

\newcommand{\der}[1]{
    \ifthenelse{\equal{\useprime}{1}}
    {#1'}
    {\dot{#1}}
}

\newcommand{\dder}[1]{
    \ifthenelse{\equal{\useprime}{1}}
    {#1''}
    {\ddot{#1}}
}

\makeatletter

\newenvironment{breakablealgorithm}
{
    \begin{center}

    \refstepcounter{algorithm}

    \hrule height .8pt depth 0pt \kern 2pt

    \renewcommand{\caption}[2][\relax]{

        {\raggedright\textbf{\ALG@name~\thealgorithm} ##2\par}

        \ifx\relax##1\relax
            \addcontentsline{loa}{algorithm}{
                \protect\numberline{\thealgorithm}##2
            }
        \else
            \addcontentsline{loa}{algorithm}{
                \protect\numberline{\thealgorithm}##1
            }
        \fi

        \kern2pt\hrule\kern2pt
    }
}
{
    \kern2pt\hrule
    \end{center}
}

\makeatother

\makeatletter
\ExplSyntaxOn
\exp_args:Nc \patchcmd { environment~Abstract~code }
  { \hbox_to_wd:nn { \z@ } { \box \g_stm_key_box } }
  { \hbox_overlap_right:n { \box \g_stm_key_box } }
  { }
  { \PackageError{paper4-layout}{CAS keyword overlay patch failed}{} }
\ExplSyntaxOff
\makeatother

\ExplSyntaxOn
\RenewDocumentCommand \printemails { }
{
  \group_begin:
  \int_compare:nNnTF { \int_use:N \g_ead_int } > { 0 }
  {
    \tex_let:D \thefootnote \relax \footnotetext
    {
      \raggedright
      \bool_if:NTF \g_stm_nologo_bool
      { 
        \int_compare:nTF { \g_ead_int = 1 }
        { \textit{Email~address:\c_space_token} }
        { \textit{Email~addresses:\c_space_token} }
      }
      { \includegraphics[height=8pt]{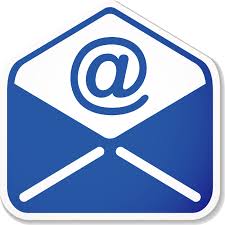}\c_space_token } 
      \seq_use:Nn \g_stm_ead_seq { ;~ }
    }
  }
  {  }
  \group_end:
}
\ExplSyntaxOff

\AddToHook{cmd/maketitle/after}{\hypersetup{pdftitle={Kernel Compensation and Quasi-Newton Iteration for Nonlinear Non-Hermitian Eigenproblems in Photonic Crystals},pdfauthor={Chenhao Jin and Hehu Xie}}}